\documentclass{article}      

\usepackage{amssymb}      
\usepackage{amsmath}      
\usepackage{amsthm}      
\usepackage{amsbsy}      
\usepackage{amscd}      
\usepackage[dvips]{graphics}

\theoremstyle{plain}      
\newtheorem{theorem}{Theorem}[section]      
\newtheorem{lemma}{Lemma}[section]      
\newtheorem{corollary}[theorem]{Corollary}      
\newtheorem{proposition}{Proposition}[section]

\newtheorem{definition}{Definition}[section]          
\theoremstyle{remark}      
\newtheorem{remark}{Remark}[section]

\newcommand{\Q}{\frac{\mathbb{Z}\left[\frac{1}{2}\right]}{\Z}}        
\newcommand{\Z}{{\mathbb{Z}}}   
   
\newcommand{\C}{{\mathbb{C}}}      
\newcommand{\R}{{\mathbb{R}}}      
\renewcommand{\H}{{\mathbb{H}}}

\newcommand{\A}{{\mathcal A}_T^{\rm ab}}

\newcommand{\al}{\widehat{\alpha}} 
\newcommand{\be}{\widehat{\beta}} 
\newcommand{\T}{T^{*}}

\newcommand{\aps}{{{\alpha^{*}}}} 
 
\newcommand{\bps}{{{\beta^{*}}}} 
 \newcommand{\als}{{\overline{\alpha}}}    
 \newcommand{\bls}{{\overline{\beta}}}

\begin{document}

\date{\today}

\title{Central extensions of the Ptolemy-Thompson group and 
quantized Teichm\"uller theory\footnote{     
This version: {June 2009}.  L.F. was partially supported by 
the ANR Repsurf:ANR-06-BLAN-0311.
      This preprint is available electronically at       
          {\rm http://www-fourier.ujf-grenoble.fr/\~{ }funar }}}      
 \author{ Louis Funar  and  Vlad Sergiescu \\      
\small \em Institut Fourier BP 74, UMR 5582 \\      
\small \em University of Grenoble I \\      
\small \em 38402 Saint-Martin-d'H\`eres cedex, France  \\      
\small \em e-mail: {\rm\{funar, sergiesc\}@fourier.ujf-grenoble.fr}  
}      

\maketitle 

\begin{abstract}
The central extension of the Thompson group  $T$ that arises in the 
quantized Teichm\"uller theory is 12 times the Euler class. 
This extension is obtained by taking a (partial) abelianization of the 
so-called braided Ptolemy-Thompson group introduced and 
studied in \cite{FK2}.  
We describe then the cyclic central extensions of $T$ by means of 
explicit presentations. 

\vspace{0.1cm}
 

\end{abstract}

\section{Introduction and statements}

Fock and Goncharov (\cite{FG1,FG}) 
improved on previous work of Faddeev, Kashaev (\cite{FaK,K}),  
Chekhov and Fock (\cite{CF}) and defined new families 
of  projective unitary representations  
of modular groupoids associated to cluster algebras. 
In the particular case of ${\rm SL}(2,\R)$ one obtained 
(projective) representations of the Ptolemy modular groupoids associated 
to triangulations of surfaces, arising in the quantification of 
the Teichm\"uller space. 
The main ingredient is the quantum dilogarithm  function which permits to deform 
the natural action of the modular groupoid on the Teichm\"uller space.
Thereby they will be called dilogarithmic representations.  
These are actually projective unitary  representations, 
or equivalently, representations of suitable central extensions.

\vspace{0.2cm}
\noindent
These representations  -- depending on a deformation parameter -- are infinite
dimensional. The general belief is that they collapse at
roots of unity to finite dimensional representations which can be identified
with the mapping class group representations arising from the quantum group
$U_q({\rm SL}(2,\R))$. Moreover, they should also coincide with the 
representations coming from the quantum hyperbolic invariants introduced 
by Baseilhac and Benedetti (see \cite{BB}).

\vspace{0.2cm}
\noindent 
There exists an universal setting for these constructions where the surface is 
the hyperbolic plane $\H^2$  (endowed with  a specific triangulation, namely 
the Farey triangulation) and the universal Teichm\"uller space 
 is the one constructed 
by Penner in \cite{pe0}. The associated modular groupoid is the 
Ptolemy groupoid  of flips on the triangulation. 
As it is well-known (see \cite{pe0}) there is a group structure 
underlying the groupoid structure that identifies the Ptolemy group of flips 
on the Farey triangulation to the 
Thompson group $T$ of  piecewise-${\rm PSL}(2,\Z)$ homeomorphisms of the circle (see \cite{CFP}). Our aim is to identify the central extension $\widehat{T}$ 
of $T$ arising in the dilogarithm representations constructed 
in (\cite{FG}, section 10, \cite{FG2}, section 3). We refer to  
$\widehat{T}$ as the dilogarithmic central extension of $T$.

\vspace{0.2cm}
\noindent 
Following \cite{GS} the cohomology ring $H^*(T)$ 
is generated by two classes $\alpha, \chi\in H^2(T)$, which are  
called the discrete Godbillon-Vey class and respectively, the Euler class. 
One can obtain $\chi$ as the Euler class of the action of 
$T$ on the circle. 

\vspace{0.2cm}
\noindent 
Moreover, central extensions of $T$ are classified up to isomorphism 
by their  extension classes in $H^2(T)$. For instance the Euler class $\chi$ 
is the extension class of the central extension 
\[ 1\to \Z \to \widetilde{T} \to T \to 1 \]
where $\widetilde{T}\subset {\rm Homeo}^+(\R)$ 
is the group of lifts of  piecewise-${\rm PSL}(2,\Z)$  homeomorphisms 
of the circle to  homeomorphisms of the real line $\R$. 

\vspace{0.2cm}
\noindent 
Our first goal is the identification of the extension class 
of $\widehat{T}$. Specifically,  
our first main result is the following: 

\begin{theorem}\label{dilog}
The class $c_{{\widehat{T}}}$ of the dilogarithmic central extension 
$\widehat{T}$
is $c_{{\widehat{T}}}= 12 \chi \in H^2(T)$. 
\end{theorem}

\vspace{0.2cm}\noindent 
Our approach consists of relating the dilogarithmic central extension 
$\widehat{T}$ to the braided Ptolemy-Thompson group $T^*$ introduced 
and studied in \cite{FK2,FK3}. This will provide a description 
of $\widehat{T}$ by means of an explicit group presentation by generators 
and relations.  The braided Ptolemy-Thompson group $T^*$ is an extension of $T$ 
by the  infinite group $B_{\infty}$ of braids, which arises as a mapping class group 
of an infinite surface.  This group is finitely presented and we are using heavily explicit 
presentations of related groups.

\vspace{0.2cm}
\noindent 
We will present a rather direct proof (using however results 
from \cite{FK2})  showing that the dilogarithmic 
extension class is a multiple of the Euler class because the  
extension  splits over the smaller Thompson group $F\subset T$. 
This multiple  is next shown to equal 12.

\vspace{0.2cm}
\noindent 
There is a general setup for studying central extensions 
of a finitely presented group in which all relations  are given arbitrary 
central lifts in the extension. In the case of the group $T$ 
this provides a series of group presentations depending on four  
integer parameters $T_{n,p,q,r}$. The dilogarithmic extension 
$\widehat{T}$ appears in this series as $T_{1,0,0,0}$. 
  
\vspace{0.2cm}\noindent 
Although our main motivation was the result of 
theorem \ref{dilog} we thought it is interesting  to obtain the  complete picture concerning 
the central extensions of $T$ in terms of explicit  group presentations. 
This amounts to understand the cohomology 2-classes of $T$ by means of their 
associated extensions.

\begin{theorem}\label{exte0}
Let $T_{n,p,q,r}$ be the group presented by  the generators  $\als, \bls, z$ and the relations: 
\[ (\bls\als)^5=z^n \]  
\[ \als^4=z^p \]
\[ \bls^3=z^q\]
\[ [\bls\als\bls, \als^2\bls\als\bls\als^2]=z^r\]
\[[\bls\als\bls, \als^2\bls\als^2\bls\als\bls\als^2\bls^2\als^2]=1\]
\[ [\als,z]=[\bls,z]=1\]
Then  each central extension of $T$ by $\Z$ is of the form $T_{n,p,q,r}$. Moreover,   
the class $c_{T_{n,p,q,r}}\in H^2(T)$ of the extension $T_{n,p,q,r}$  
is given by:  
\[ c_{T_{n,p,q,r}}= (12n-15p-20q-60r) \chi + r \alpha \]
\end{theorem}	

\vspace{0.2cm}\noindent     
The extension classes behave linearly on the parameters and, in order
to find their coefficients, one has to  
make use of several explicit (central) extensions. 
There are only a few such extensions 
in the literature, namely those coming  from mapping class groups of 
planar surfaces with infinitely many punctures. For general central extensions
we face new complications due to the presence of the discrete Godbillon-Vey class. 
The Ptolemy-Thompson group 
$T^*$ is one such mapping class group, but it is useless 
in computing the coefficients of the Godbillon-Vey class. 
The other extension of the same kind 
is the Greenberg-Sergiescu acyclic extension (see \cite{GrS}) and hence 
the core of the proof   of theorem \ref{exte0} 
consists of computations within the Greenberg-Sergiescu extension. 
This lead us to the formula for $c_{T_{n,p,q,r,0}}$   stated in theorem \ref{exte0}.

\vspace{0.2cm}\noindent  
All over this paper the Ptolemy-Thompson group $T$ appears in every one of its 
instances, as a group of homeomorphisms of the circle, as the group of flips, 
as the group $PPSL(2,\Z)$ and as a mapping class group. 
A main difficulty is to pass from one description to another one. 

\vspace{0.2cm}\noindent 
The plan of the paper is as follows. 
We introduce in the first part  the 
braided Ptolemy-Thompson  $T^*$ and the  induced 
abelianized central extension $T^*_{\rm ab}$. We will explain that the class 
$c_{T^*_{\rm ab}}$ is a multiple of the Euler class and provide a short  
proof that it should be 12 times the later. In particular we describe 
all central extensions associated to multiples of the Euler class. 
In the second part we give a quick overview of the quantization of the Teichm\"uller space,  
following mostly Fock and Goncharov's series of papers, to  introduce the dilogarithmic 
projective representation of the Ptolemy-Thompson group.  
As projective representations are in one-to-one correspondence with 
central extensions one obtains what we will call the dilogarithmic extension $\widehat{T}$ 
of $T$. In the final section of this part we prove that  
$\widehat{T}$ and  $T^*_{\rm ab}$ are isomorphic. 
The third part of the paper is devoted to the classification of all central extensions 
of $T$ by $\Z$. The major problem is to understand the behavior of the 
discrete Godbillon-Vey class in extensions. 
If the braided Ptolemy-Thompson group $T^*$ was the geometric model 
in the first part now we make use of the Greenberg-Sergiescu ${\mathcal A}_T$ of $T$.  
The method for computing the classes of extensions   
is inspired by Milnor's algorithm concerning the Euler class 
of a surface group representation. 
The last section contains some speculations on what should be 
the geometric extensions  of $T$ on one hand, and 
the central extensions of  the mapping class group of 
a finite (punctured) surface arising in the quantized 
Teichm\"uller theory. 

\vspace{0.2cm}
\noindent 
Although we use some of the results from \cite{FK2} and \cite{GrS}, quite a 
large amount of these papers is not essential for understanding the present paper.
We  also outlined the main constructions from \cite{FG1} in the case of the 
universal Teichm\"uller space in order to make the paper more self-contained.

\vspace{0.2cm}
\noindent 
{\bf Acknowledgements.} The authors are indebted to St\'ephane Baseilhac, 
Volodya Fock, Rinat Kashaev, 
Christophe Kapoudjian and Greg McShane for useful discussions.

\tableofcontents

\section{The central extension induced from the braided Ptolemy-Thompson group 
$T^*$}
\subsection{The Thompson groups $F$ and $T$}
The Thompson group $F$ is the group of dyadic piecewise affine 
homeomorphism of $[0,1]$. We describe every element $\gamma$ of $F$ 
as follows. There exist two partitions 
of $[0,1]$ into consecutive intervals $I_1,I_2,\ldots, I_k$ on one side 
and $J_1,J_2,\ldots, J_k$ whose end points are dyadic numbers and 
$\gamma$ sends affinely each interval $I_j$ into its counterpart $J_j$, for 
all $j\in \{1,2,\ldots, k\}$. Thus the restriction of $\gamma$ to any 
interval $I_j$ is given by $\gamma(x)=a_jx+b_j$, where 
$a_j=2^{n_j}, n_j\in \Z$ and $b_j$ belong to the set of dyadic numbers 
i.e. $b_j=\frac{p_j}{2^{m_j}}$, $p_j,m_j\in \Z$. 

\vspace{0.2cm}\noindent 
The group $F$ is generated by the two 
elements $A$ and $B$ below: 
\[ A(x)=\left\{ \begin{array}{lll}
\frac{x}{2}, & \mbox{ if }  & x\in[0,\frac{1}{2}] \\
x-\frac{1}{4}, & \mbox{ if }  & x\in[\frac{1}{2}, \frac{3}{4}] \\
2x-1, & \mbox{ if }  & x\in[\frac{3}{4}, 1] 
\end{array}\right., \,\,\, 
 B(x)=\left\{ \begin{array}{lll}
x, & \mbox{ if }  & x\in[0,\frac{1}{2}] \\
\frac{x}{2}+\frac{1}{4}, & \mbox{ if }  & x\in[\frac{1}{2}, \frac{3}{4}] \\
x-\frac{1}{8}, & \mbox{ if }  & x\in[\frac{3}{4}, \frac{7}{8}] \\
2x-1,\mbox{ if }  & x\in[\frac{7}{8}, 1] 
\end{array}\right. 
\] 

\vspace{0.2cm}\noindent 
Moreover, the group $F$ has the presentation  
\[ F=\langle A, B; [AB^{-1}, A^{-1}B A]=1, \,\,  [AB^{-1}, A^{-2}B A^2]=1\rangle \]

\vspace{0.2cm}\noindent 
There is a geometric encoding of elements of $F$ as pairs of 
(stable) finite rooted binary trees. Each finite rooted binary tree 
encodes a subdivision of $[0,1]$ into dyadic intervals. Adding two descending 
edges to a vertex amounts to subdivide the respective interval into two equal 
halves. Given the two subdivisions the element of $F$ sending 
one into the other is uniquely determined. The pair of trees is determined 
up to stabilization, namely adding extra couples of descending edges to 
corresponding vertices in both trees. 
 
\begin{center}
\includegraphics{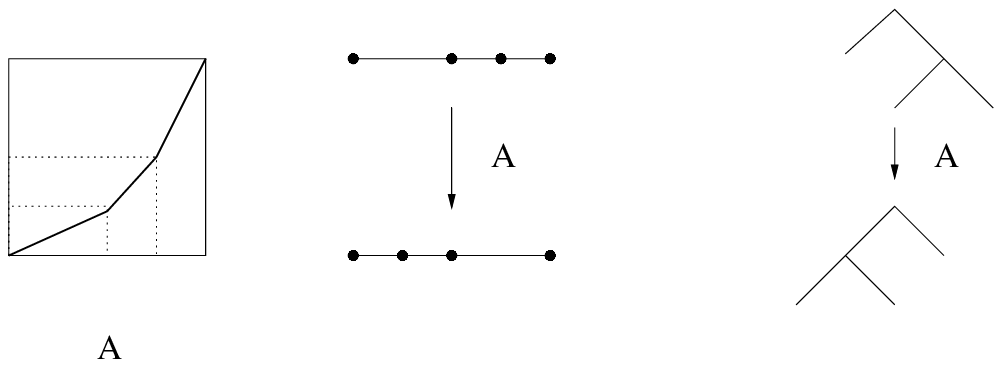}
\end{center}
\begin{center}
\includegraphics{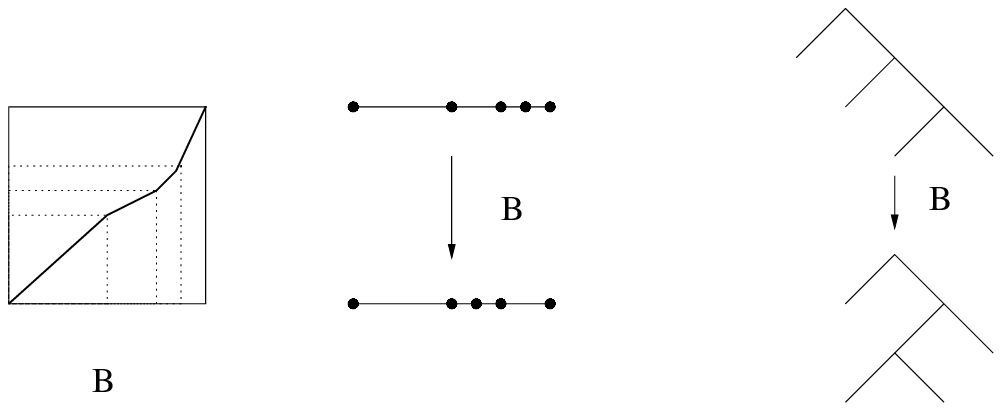}
\end{center}

\vspace{0.2cm}\noindent 
The Thompson group $T$ is the group of dyadic piecewise affine 
homeomorphisms of $S^1=[0,1]/0\sim 1$. It is not hard to see that 
$A,B$ and $C$ generate $T$, where 
\[
C(x)=\left\{ \begin{array}{lll}
\frac{x}{2}+\frac{3}{4}, & \mbox{ if }  & x\in[0,\frac{1}{2}] \\
2x-1, & \mbox{ if }  & x\in[\frac{1}{2}, \frac{3}{4}] \\
x-\frac{1}{4}, & \mbox{ if }  & x\in[\frac{3}{4}, 1] 
\end{array}\right. 
\] 

\vspace{0.2cm}\noindent 
The presentation of $T$ in terms of the generators 
$A,B,C$ consists  of the two relations above with  four more 
relations to be added: 
\[ T=\left\langle \begin{array}{cl} 
A,B,C; &  CA=(A^{-1}CB)^2,\,\; 
(A^{-1}CB)(A^{-1}BA)=B(A^{-2}CB^2), \\
 &C^3=1, \;\, C=B A^{-1}CB,\ [AB^{-1}, A^{-1}B A]=1, \,\;\; [AB^{-1}, A^{-2}B A^2]=1
\end{array}\right
\rangle \]

\vspace{0.2cm}\noindent 
We can associate a pair of  (stable) trees to encode an element of $T$, as well, but we have to specify additionally, where the origin is sent to. 
We denote the image of the origin in the trees by marking the 
leftmost leaf of its domain. Usually the origin is fixed to be the
leftmost leaf in the first tree. For instance, the element $C$ 
has the following description:

\begin{center}
\includegraphics{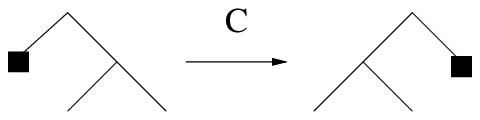}
\end{center}

\begin{remark}\label{presT}
The group $T$ has also another presentation with generators 
$\alpha$ and $\beta$ and relations 
\[ \alpha^4=\beta^3=1\]
\[ [\beta\alpha\beta, \alpha^2\beta\alpha\beta\alpha^2]=1\]
\[  [\beta\alpha\beta, \alpha^2\beta^2\alpha^2\beta\alpha\beta\alpha^2
\beta\alpha^2]=1\]
\[ (\beta\alpha)^5=1\]
\vspace{0.2cm}\noindent 
If we set $A=\beta\alpha^2$, $B=\beta^2\alpha$ and 
$C=\beta^2$ then we obtain the generators $A,B,C$ of the group $T$, 
considered above. Then the  two commutativity relations above 
are equivalent to 
\[ [AB^{-1}, A^{-1}B A]=1, \,\;\; [AB^{-1}, A^{-2}B A^2]=1\]
\end{remark}

\subsection{Mapping class groups of infinite surfaces}
\begin{definition}
A {\em rigid structure} $d$  on the surface 
$\Sigma$ is a decomposition  of $\Sigma$ into 2-disks with disjoint 
interiors, called elementary pieces. We suppose that the closures 
of the elementary pieces are still 2-disks.  

\vspace{0.2cm}
\noindent
We assume that we are given a family $F$, called the family of admissible sub-surfaces of $\Sigma$, of compact sub-surfaces of $\Sigma$ 
such that each member of $F$ is a finite union of elementary pieces. 
\end{definition}

\vspace{0.2cm}
\noindent
Given the data $(\Sigma, d, F)$ we can associate the asymptotic 
mapping class group ${\cal M}(\Sigma, d, F)$ as follows. We restrict first to 
those homeomorphisms that act in the simplest possible way at infinity. 

\begin{definition}
A homeomorphism $\varphi$ between two surfaces endowed 
with rigid structures is {\em rigid} if it sends the rigid structure of one 
surface onto the rigid structure of the other. 

\vspace{0.2cm}
\noindent
The homeomorphism $\varphi:\Sigma\to \Sigma$ is  said to be 
{\em asymptotically rigid}  
if  there exists some admissible subsurface 
$C\subset \Sigma$, called a support for $\varphi$, 
such that $\varphi(C)\subset \Sigma$ is 
also an admissible subsurface of $\Sigma$ and  
the restriction $\varphi|_{\Sigma-C}:\Sigma-C\to \Sigma-\varphi(C)$ 
is rigid. 
\end{definition}

\vspace{0.2cm}
\noindent As it is customary when studying mapping class groups  
we consider now isotopy classes of such homeomorphisms. 

\begin{definition}
The group ${\cal M}(\Sigma, d, F)$ of isotopy classes of asymptotically rigid homeomorphisms is  called the 
{\em asymptotic mapping class group} of $\Sigma$ corresponding to the  
rigid structure $d$ and family  of admissible sub-surfaces $F$. 
\end{definition}

\begin{remark}
Two asymptotically rigid homeomorphisms that are isotopic should be 
isotopic among asymptotically rigid homeomorphisms. 
\end{remark}

\subsection{$T$ and $T^*$ as mapping class groups}

The surfaces below will be oriented and all      
homeomorphisms considered in the sequel will be      
orientation-preserving, unless the opposite is explicitly stated. 
Actions in the sequel are left actions and the composition of maps 
is the usual one, namely we start composing from right to the left.    
\begin{definition}\label{ss}      
The ribbon tree $D$ is the planar surface obtained by   
thickening in the plane the infinite binary        
tree. We denote by $D^{*}$ the ribbon tree  
with infinitely many punctures, one puncture for each edge of the 
tree.  A homeomorphism of $D^*$ is a homeomorphism of $D$ which permutes the punctures $D^*$.
\end{definition}

\begin{definition}
A {\em rigid structure} on $D$ or $D^{*}$ is a decomposition  into
hexagons by means of a family of arcs whose 
endpoints are on the boundary of $D$. Each hexagon contains no puncture within its interior but each 
 arc passes through a unique puncture in the case of $D^{*}$. It is assumed that these arcs  
are pairwise non-homotopic in $D$, by homotopies keeping the boundary points of the arcs 
on the boundary of $D$. The choice of a rigid structure of reference is called the {\em canonical rigid structure}. 
The canonical rigid structure 
of the ribbon tree $D$ is such that each arc of this rigid structure crosses once and transversely a unique edge of the tree.
 The canonical 
rigid structure on $D^{*}$ is assumed to coincide with the canonical rigid structure of $D$ when forgetting the punctures. 
See the Figure below: 

  \begin{center}   
\includegraphics{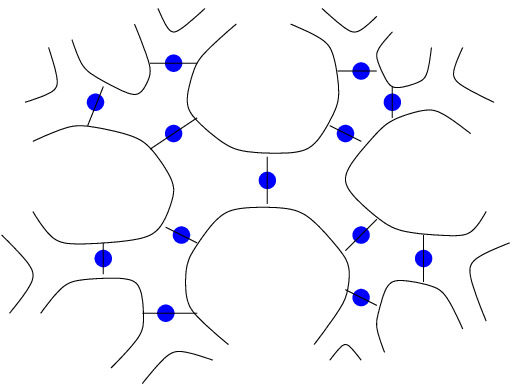}      
\end{center}  
\noindent     
A planar subsurface of $D$ or $D^*$ is {\em admissible}  if  
it is a connected finite union of hexagons coming from the canonical rigid 
structure.  

\vspace{0.2cm}\noindent 
One denotes by $T$ and $T^{*}$ the group of isotopy classes of asymptotically rigid  homeomorphisms of $D$ and $D^{*}$, respectively.
\end{definition}

\begin{remark} 
There exists a cyclic order on the frontier arcs of an  
admissible subsurface induced by the planarity. An asymptotically 
rigid homeomorphism necessarily preserves the cyclic order  
of the frontier for any admissible subsurface.  

\end{remark}

\vspace{0.2cm}\noindent 
The mapping class group $T$ is isomorphic to the Thompson group which is commonly denoted $T$. This fact has been widely developed in \cite{KS} and \cite{FK2}.
We consider the following elements of $T$, defined as mapping classes 
of asymptotically rigid homeomorphisms: 
\begin{itemize} 
\item The support of the element $\beta$ is the central hexagon on the figure below.
Further,  $\beta$ acts as the counterclockwise rotation of order  
three which permutes the three branches 
of the ribbon tree issued from the hexagon.  

\begin{center} 
\includegraphics{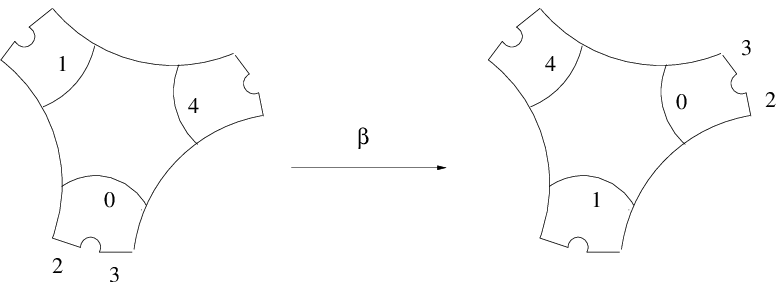}      
\end{center}

In fact, $\beta$ is globally rigid.

\item  The support of $\alpha$  is the union of two adjacent hexagons, 
one of them being the support of $\beta$.  Then $\alpha$  
rotates counterclockwise the support of angle $\frac{\pi}{2}$, by  
permuting the four branches of the ribbon tree
issued from the support.  

\begin{center} 
\includegraphics{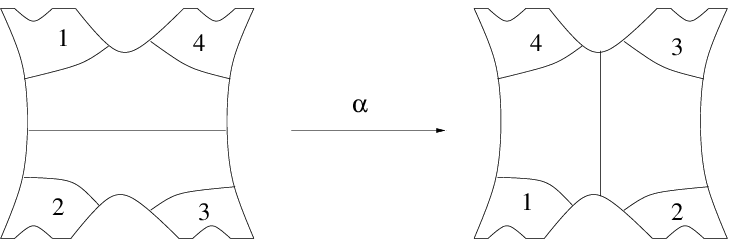}      
\end{center}

Note that $\alpha$ is not globally rigid, but $\alpha^2$ is.
\end{itemize}

\subsection{The relative abelianization of the 
braided Ptolemy-Thompson group $T^*$}
Recall from \cite{FK2,FK3} that there exists a natural 
surjection homomorphism $T^*\to T$ between the two mapping class groups, 
which is obtained by {\em forgetting the punctures}.  
Its kernel is the infinite braid group $B_{\infty}$ 
consisting of those braids in the punctures 
of $D^*$ that move non-trivially only finitely many punctures. 
In other words $B_{\infty}$ is the direct limit of an ascending 
sequence of braid groups associated to an exhaustion of $D^*$ by punctured 
disks.  This yields the  following exact sequence description of $T^*$:
\[ 1 \to B_{\infty} \to \T\to T \to 1 \]

\vspace{0.2cm}
\noindent 
Observe that $H_1(B_{\infty})=\Z$. Thus, the abelianization homomorphism  
$B_{\infty}\to H_1(B_{\infty})=\Z$  
induces a central extension $T^*_{\rm ab}$ of $T$, where one 
replaces $B_{\infty}$ by its abelianization 
$H_1(B_{\infty})$, as in the diagram below: 

\[ \begin{array}{ccccccc}
1  \to & B_{\infty} & \to & \T         & \to &  T        & \to 1 \\
       & \downarrow &     & \downarrow &     & \parallel &   \\ 
1  \to & \Z         & \to & T^*_{\rm ab}          & \to &   T       & \to 1  
\end{array}
\]

\vspace{0.2cm}
\noindent Then $T^*_{\rm ab}$ is the relative abelianization of $T^*$ over $T$. 
We are not only able to make computations in the 
mapping class group $T^*$ and thus in $T^*_{\rm ab}$,  but also 
to interpret the algebraic relations in $T^*_{\rm ab}$ in geometric terms. 
 
\begin{proposition}\label{abelcar}
The group $T^*_{\rm ab}$ has the presentation with three generators 
$\alpha^*_{\rm ab}$, $\beta^*_{\rm ab}$ and $z$ and the relations 
\[ {\alpha^*_{\rm ab}}^4={\beta^*_{\rm ab}}^3=1, (\beta^*_{\rm ab}\alpha^*_{\rm ab})^5=z, 
[\alpha^*_{\rm ab}, z]=1, [\beta^*_{\rm ab}, z]=1\]
\[ [\beta^*_{\rm ab}\alpha^*_{\rm ab}\beta^*_{\rm ab} \, , \, {\alpha^*_{\rm ab}}^2\beta^*_{\rm ab}\alpha^*_{\rm ab}\beta^*_{\rm ab}{\alpha^*_{\rm ab}}^2]=
[\beta^*_{\rm ab}\alpha^*_{\rm ab}\beta^*_{\rm ab} \, , \, {\alpha^*_{\rm ab}}^2{\beta^*_{\rm ab}}^2 
{\alpha^*_{\rm ab}}^2\beta^*_{\rm ab}\alpha^*_{\rm ab}\beta^*_{\rm ab}{\alpha^*_{\rm ab}}^2{\beta^*_{\rm ab}}
{\alpha^*_{\rm ab}}^2]=1\]
Moreover the projection map $T^*_{\rm ab}\to T$ sends $\alpha^*_{\rm ab}$ to $\alpha$, 
$\beta^*_{\rm ab}$ to $\beta$ and $z$ to identity.  
\end{proposition}
\begin{proof}
Recall from \cite{FK2} that $T^*$ is generated by two elements $\alpha^*$ and 
$\beta^*$ below. 
\begin{itemize} 
\item The support of the element $\beta^{*}$ of $T^{*}$ is the central hexagon.  Further $\beta^{*}$ acts as the counterclockwise rotation of order  
three which permutes cyclically  
the punctures.  One has ${\beta^{*}}^{3}=1$.

\begin{center}
\includegraphics{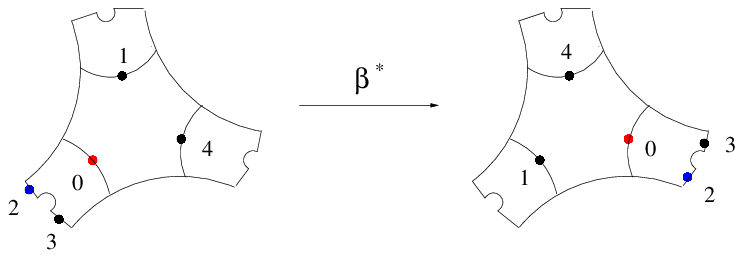}
\end{center}

\item  The support of the element $\alpha^{*}$ of $T^{*}$ is the union 
of two adjacent hexagons, 
one of them being the support of $\beta^{*}$.  
Then $\alpha^{*}$ rotates counterclockwise the support 
of angle $\frac{\pi}{2}$, by  
keeping fixed the central puncture. One has  ${\alpha^{*}}^4=1$. 

\begin{center}
\includegraphics{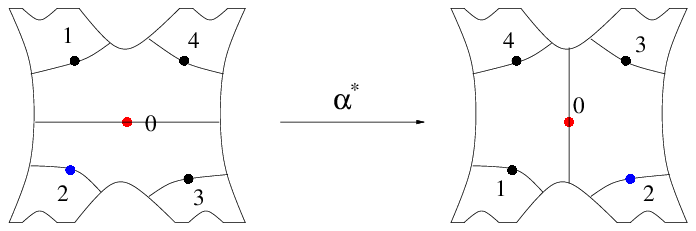}
\end{center}
\end{itemize}

\vspace{0.2cm}
\noindent
Let now $e$ be a simple arc in  $D^{*}$ which connects two punctures. 
We associate a braiding $\sigma_e\in B_{\infty}$ to $e$ by 
considering the homeomorphism
that moves clockwise the punctures at the 
endpoints of the edge $e$ in a small neighborhood of the edge,  
in order to interchange their positions. This means that, if $\gamma$ is an arc transverse to $e$, then the braiding $\sigma_{e}$
moves $\gamma$ on the left when it approaches $e$. Such a braiding will be called {\it positive}, while $\sigma_{e}^{-1}$ is
{\it negative}.

\begin{center} 
\includegraphics{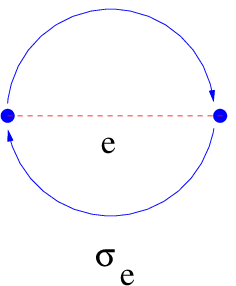}      
\end{center}

\vspace{0.2cm}
\noindent
It is known that $B_{\infty}$ is generated by the braids $\sigma_e$ 
where $e$ runs over the edges of the binary tree with vertices 
at punctures. Let $\iota:B_{\infty}\to T^*$ be the inclusion. 
It is proved in \cite{FK2} that 
the braid generator $\sigma_{[02]}$ associated to the edge joining the 
punctures numbered $0$ and $2$ has the image 
\[ \iota(\sigma_{[02]})=(\bps\aps)^5\]
because we have 

\begin{center}
\includegraphics{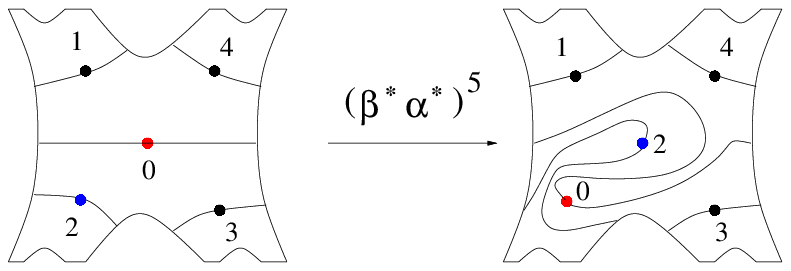}
\end{center}

\vspace{0.2cm}
\noindent
Recall next that all braid generators $\sigma_e$ are conjugate and 
call $z$ their image in $T^*_{\rm ab}$. It follows that $T^*_{\rm ab}$ is an 
extension of $T$ by $\Z$. Moreover, it is simple to check that 
$\aps\sigma_{[02]}\alpha^{-1}$ is also a braid generator, namely 
$\sigma_{[\aps(0)\aps(2)]}$. The same holds true for 
$\bps\sigma_{[02]}\bps^{-1}=\sigma_{[\bps(0)\bps(2)]}$. 
This implies that the extension $T^*_{\rm ab}$ is central.  

\vspace{0.2cm}
\noindent
In particular, a presentation of $T^*_{\rm ab}$ can be obtained by  looking at 
the lifts of  relations in $T$, together with those coming from 
the fact that $z$ is central. 

\vspace{0.2cm}
\noindent
The first  set relations above are obviously satisfied by $T^*_{\rm ab}$. 
Finally recall from \cite{FK2} that $T^*$ splits over the smaller 
Thompson group $F$ and thus the following relations hold true in $T^*$:
\[ [\beta^*\alpha^*\beta^* \, , \, {\alpha^*}^2\beta^*\alpha^*\beta^*{\alpha^*}^2]=
[\beta^*\alpha^*\beta^* \, , \, {\alpha^*}^2{\beta^*}^2 
{\alpha^*}^2\beta^*\alpha^*\beta^*{\alpha^*}^2{\beta^*}
{\alpha^*}^2]=1\]
Thus the second set of relations are  
automatically verified in $T^*_{\rm ab}$. 
Since these relations form a complete set of lifts of  relations 
presenting $T$ and $z$ is central, then they represent a complete 
system of relations in $T^*_{\rm ab}$. This ends the proof. 
\end{proof}

\subsection{Computing the class of $T^*_{\rm ab}$}

\begin{lemma}\label{multiple}
The class $c_{T^*_{\rm ab}}$ is a multiple of the Euler class. 
\end{lemma}
\begin{proof}
Since $T^*$ splits over Thompson group $F\subset T$ (see \cite{FK2}) it follows that 
$T^*_{\rm ab}$ also splits over $F$. Therefore the extension class  $c_{T^*_{\rm ab}}$ 
lies in the kernel of the restriction map $H^2(T)\to H^2(F)$. According to \cite{GS} 
the kernel is generated by the Euler class. 
\end{proof}

\vspace{0.2cm}
\noindent
Let us introduce the group $T_{n,p,q}$ presented by  the generators  
$\als, \bls, z$ and the relations: 
\[ (\bls\als)^5=z^n \]  
\[ \als^4=z^p \]
\[ \bls^3=z^q\]
\[ [\bls\als\bls, \als^2\bls\als\bls\als^2]=1\]
\[[\bls\als\bls, \als^2\bls\als^2\bls\als\bls\als^2\bls^2\als^2]=1\]
\[ [\als,z]=[\bls,z]=1\]
Recall from proposition \ref{abelcar} 
that $T^*_{\rm ab}=T_{1,0,0}$. It is easy to see that $T_{n,p,q}$ are central extensions 
of $T$ by $\Z$.  Because of the last two commutation relations the extension 
$T_{n,p,q}$ splits over the Thompson group $F$. 
Thus the restriction of $c_{T_{n,p,q}}$ to $F$ vanishes and 
a fortiori the restriction to the commutator subgroup  $F'\subset F$.   
According to \cite{GS} we have $H^2(F')=\Z\alpha$ where 
$\alpha$ is the discrete  Godbillon-Vey class. Thus the map 
$H^2(T)\to H^2(F')$ is the projection $\Z\alpha\oplus \Z\chi\to \Z\alpha$. 
Since $c_{\widehat{T}}$ belongs to the 
kernel of $H^2(T)\to H^2(F')$. This proves that 
$c_{T_{n,p,q}}\in \Z\chi$. Set $c_{T_{n,p,q}}=\chi(n,p,q) \chi$.

\begin{proposition}\label{linearform}
We have $\chi(n,p,q)= 12n-15p-20q$. 
\end{proposition}

\begin{proof} 
Observe first that:  
\begin{lemma}\label{linearm}
The function $\chi(n,p,q):\Z^3\to \Z$ is linear. 
\end{lemma}
\begin{proof}
The group $T$ contains all finite cyclic subgroups. Let us fix some 
positive integer $k$ and set $y_k\in T$ for an element of order $k$. 
Denote by $G_k\subset T$ the cyclic subgroup generated by $y_k$. 

\vspace{0.2cm}
\noindent
We want to compute the restriction $c_{T_{n,p,q}}|_{G_k}$, 
namely the image of $c_{T_{n,p,q}}$ under the obvious 
morphism $H^2(T)\to H^2(G_k)$. 
Since $G_k$ is cyclic $H^2(G_k)$ is the cyclic $\Z/k\Z$ of order $k$, generated 
by the restriction of the Euler class $\chi|_{G_k}$. 
In particular,  
\[ c_{T_{n,p,q}}|_{G_k}=\chi(n,p,q)\chi|_{G_k} \in \Z/k\Z= H^2(G_k) \]
On the other hand we can compute  the 
Euler class of a central extension of a cyclic subgroup 
of $T\subset {\rm Homeo}^+(S^1)$  by means of the following Milnor-Wood 
algorithm. Assume that $y_k=w_k(\alpha,\beta)$ is given by a word 
in the generators $\alpha,\beta$ of $T$. Consider next 
the word $\overline{y_k}=w_k(\als,\bls)\in T_{n,p,q}$. 
Then $\overline{y_k}^k=z^{e}$ for some $e$ and the Euler class 
is given by the value of $e$ modulo $k$. 

\vspace{0.2cm}
\noindent However  $\overline{y_k}^k$ is an element of the center of 
$T_{n,p,q}$ and thus it is a product of conjugates of 
$(\bls\als)^5$, $\als^4$, $\bls^3$ and the commutation 
relations (which are not involving $z$). If $\varphi_{n,p,q}(w)$  is the additive map that 
associates to any sub-word of $w$  from the subset  
$\{(\bls\als)^5, \als^4, \bls^3\}$ respectively $n, p$ and $q$ 
it follows that 
\[ \varphi_{n,p,q}(w)=n\varphi_{1,0,0}(w)+p\varphi_{0,1,0}(w)+
q\varphi_{0,0,1}(w) \]
In particular, by choosing $w=w_k(\als,\bls)^k$ we obtain that 
the coefficient of the Euler class $\chi(n,p,q,r)$ verifies 
\[ \chi(n,p,q)=n\chi(1,0,0)+p \chi(0,1,0)+
q \chi(0,0,1) \,\, ({\rm mod } \, k) \]
This equality holds for all  natural $k$ and therefore $\chi(n,p,q)$ is linear. 
\end{proof}

\vspace{0.2cm}
\noindent
{\em End of the proof of the proposition \ref{linearform}}. 
By the definition of the extension class $\chi(n,p,q)$ takes the same 
value for those $T_{n,p,q}$ that are isomorphic by an isomorphism 
inducing identity on $T$ and on the center. 
Such isomorphisms $L$ have the form $L(\als)=\als z^x$
and $L(\bls)=\bls z^y$. Thus $T_{n,p,q}$ is isomorphic as extension 
to $T_{n',p',q'}$ if and only if there exists integers $x,y\in Z$ such that 
\[ n=n'+5x+5y, \, p=p'+4x, \, q=q'+3y \]
In particular the linear form $\chi(n,p,q)$ should be invariant by the  
transforms corresponding to arbitrary $x,y\in \Z$ and thus it should be a multiple of $12n-15p-20q$. 

\vspace{0.2cm}
\noindent Last, the extension $\widetilde{T}\subset{\rm Homeo}^+(\R)$ 
is known to have Euler number 1. The central element $z$ acts as 
the unit translation on the line and in order to identify it with an element of 
the family $T_{n,p,q}$ it suffices to compute the 
rotation index associated to the elements $\alpha^4,\beta^3$ and 
$(\beta\alpha)^5$. One obtains that $\widetilde{T}$ is actually isomorphic 
to $T_{3,1,1}$. Thus $\chi(n,p,q)$ is precisely given by 
the claimed linear form. 
\end{proof}

\begin{corollary}\label{12}
We have $c_{T^*_{\rm ab}}=12\chi$. 
\end{corollary}

\begin{remark}
The extension $T^*\to T$ splits also 
over the subgroup $\langle \alpha^2,\beta\rangle$ which is isomorphic to $PSL(2,\Z)$ 
(see \cite{FK2}). This implies that $c_{T^*_{\rm ab}}$ is a multiple of $6\chi$.  
\end{remark}

\section{Quantum Teichm\"uller space and dilogarithmic representations}
\subsection{The Ptolemy-Thompson group and the Ptolemy groupoid}
We will use the terms triangulation of $\H^2$ for 
ideal  {\em locally finite} triangulations of the hyperbolic space $\H^2$ i.e. 
a countable locally finite collection of geodesics whose 
complementary regions are triangles. The {\em vertices} of the triangulation 
are the asymptotes of the constitutive geodesics, viewed as points of the 
circle at infinity $S^1_{\infty}$. 

\vspace{0.2cm}
\noindent   
Consider the  basic ideal 
triangle having vertices at $1,-1, \sqrt{-1}\in S^1_{\infty}$ in the 
unit disk model  $\mathbb D$ of $\H^2$. The orbits of its sides 
by the group $PSL(2,\Z)$ is the so-called {\em Farey triangulation} $\tau_0$.

\vspace{0.2cm}
\noindent  
A triangulation $\tau$ is {\em marked} if one fixes  
a distinguished oriented edge (abbreviated d.o.e.) $\vec{a}$ of it.  
The standard marking of the Farey triangulation $\tau_0$ is the oriented edge 
$\vec{a_0}$ joining $-1$ to $1$.

\vspace{0.2cm}
\noindent 
We define next a {\em marked tessellation} of $\H^2$ to be an equivalence class 
of marked triangulations of $\H^2$ with respect to the action of $PSL(2,\R)$. 
Since the action of $PSL(2,\R)$ is 3-transitive each tessellation can be 
uniquely represented by its {\em canonical marked triangulation} 
containing the basic ideal triangle and whose d.o.e. is $\vec{a_0}$. 
The marked tessellation is of Farey-type if its canonical marked triangulation 
has the same vertices  and all but finitely many triangles (or sides) as the Farey triangulation. 
Unless explicitly stated otherwise all tessellations considered 
in the sequel will be Farey-type tessellations.  In particular, 
the ideal triangulations have the same vertices as $\tau_0$ and coincide 
with $\tau_0$ for all but finitely many ideal triangles.

\begin{definition}
The objects of the {\em (universal) Ptolemy groupoid} $Pt$ are the 
marked tessellations of Farey-type. The morphisms  between  two objects 
$(\tau_1,\vec{a_1})$ and $(\tau_2,\vec{a_2})$ are {\em eventually trivial} 
permutations maps $\phi:\tau_1\to \tau_2$ such that 
$\phi(\vec{a_1})=\vec{a_2}$. When marked tessellations are represented 
by their canonical triangulations, the later coincide for all 
but finitely many triangles. Recall that $\phi$ is said to be eventually trivial 
if the induced correspondence at the level of canonical triangulations
is the identity for all but finitely many edges. 
\end{definition}

\vspace{0.2cm}
\noindent 
We define now particular elements of $Pt$ called flips. 
Let $e$ be an edge of the marked tessellation represented by 
the canonical marked triangulation  $(\tau, \vec{a})$.  
The result of the flip $F_e$ on $\tau$ is the triangulation
$F_e(\tau)$ obtained  from $\tau$ by changing only 
the two neighboring triangles containing the edge $e$, according 
to the picture below: 

\vspace{0.2cm}
\begin{center}
\includegraphics{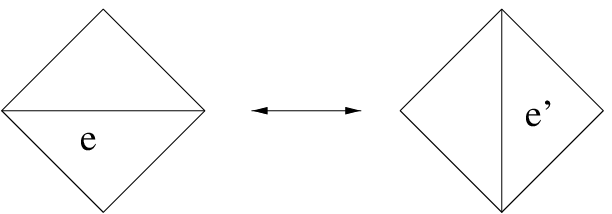}
\end{center}

\vspace{0.2cm}
\noindent 
This means that we remove $e$ from $\tau$ and further add the 
new edge $e'$ in order to get $F_e(\tau)$.  
In particular there is a natural correspondence $\phi:\tau\to F_e(\tau)$ 
sending $e$ to $e'$ and identity for all other edges. The result of a flip 
is the new triangulation together with this edges' correspondence.

\vspace{0.2cm}
\noindent 
If $e$ is  not the d.o.e. of $\tau$ then $F_e(\vec{a})=\vec{a}$. 
If $e$ is the d.o.e. of $\tau$ then $F_e(\vec{a})=\vec{e'}$, where 
the orientation of $\vec{e'}$ is chosen so that 
the frame $(\vec{e}, \vec{e'})$ is positively oriented.

\vspace{0.2cm}
\noindent We define now the flipped tessellation $F_e((\tau, \vec{a}))$ 
to be the tessellation $(F_e(\tau), F_e(\vec{a}))$. 
It is proved in \cite{pe0} that flips generate the Ptolemy groupoid i.e. 
any element of $Pt$ is a composition of flips.

\vspace{0.2cm}
\noindent 
For any marked tessellation $(\tau, \vec{a})$ there is defined a 
characteristic map $Q_{\tau}:{\mathbb Q}-\{-1,1\}\to \tau$. 
Assume that $\tau$ is the canonical triangulation representing this 
tessellation. 
We label first by $\mathbb Q\cup{\infty}$ 
the vertices of $\tau$, by induction:
\begin{enumerate}
\item $-1$ is labeled by $0/1$, $1$ is labeled by $\infty=1/0$ and 
$\sqrt{-1}$ is labeled by $-1/1$. 
\item If we have a triangle in $\tau$ having two vertices already 
labeled by $a/b$ and $c/d$  then its third vertex is labeled $(a+c)/(b+d)$. 
Notice that vertices in the upper half-plane are labeled by negative 
rationals and those from the lower half-plane by positive rationals. 
\end{enumerate}
\vspace{0.2cm}
\noindent 
As it is well-known this labeling produces a bijection between 
the set of vertices of $\tau$ and $\mathbb Q\cup{\infty}$. 

\vspace{0.2cm}
\noindent 
Let now  $e$ be an edge of $\tau$, which is different from $\vec{a}$.
Let $v(e)$ be the vertex  opposite to $e$ of the triangle $\Delta$ of $\tau$ 
containing $e$ in its frontier and lying in the component of 
${\mathbb D}-e$ which does not contain 
$\vec{a}$. We associate then to $e$ the label of $v(e)$. 
Give also $\vec{a}$ the label $0\in {\mathbb Q}$.  
This way one obtains a bijection  
$Q_{\tau}:{\mathbb Q}-\{-1,1\}\to \tau$. 

\vspace{0.2cm}
\noindent 
Remark that if $(\tau_1,\vec{a_1})$ and $(\tau_2,\vec{a_2})$ are 
marked tessellations then there exists a unique map $f$ between their vertices 
sending triangles to triangles and marking on the marking. 
Then $f\circ Q_{\tau_1}=Q_{\tau_2}$. 

\vspace{0.2cm}
\noindent 
The role played by $Q_{\tau}$ is to allow flips to be indexed by 
the rationals and not on edges of $\tau$. 

\begin{definition}
Let ${\rm FT}$ be the set of marked tessellations of Farey-type. 
Define the action of the free monoid  $M$ generated by ${\mathbb Q}-\{-1,1\}$
on ${\rm FT}$ by means of:  
\[ q \cdot (\tau, \vec{a}) = F_{Q_{\tau}(q)}(\tau, \vec{a}),\,  \mbox{ for } \, 
q\in {\mathbb Q}-\{-1,1\}, (\tau, \vec{a})\in {\rm FT}\]
We set  $f\sim f'$ on $M$  if the  two actions of $f$ and $f'$ on  
${\rm FT}$ coincide. Then the induced composition law 
on  $M/\sim$ is a monoid structure for which each element has an inverse. This  
makes  $M/\sim$ a group, which is called the Ptolemy group $T$ (see \cite{pe0} for more details).  
\end{definition}

\vspace{0.2cm}
\noindent 
In particular it makes sense to speak of flips in the present case
and thus flips generate the  
Ptolemy group. 

\vspace{0.2cm}
\noindent 
The notation $T$ for the Ptolemy group is not misleading because 
this group is isomorphic to the Thompson group $T$ and for this reason, 
we preferred to call it the Ptolemy-Thompson group.

\vspace{0.2cm}
\noindent 
Given two marked tessellations $(\tau_1,\vec{a_1})$ and $(\tau_2,\vec{a_2})$ 
the  combinatorial isomorphism $f:\tau_1\to \tau_2$ from above  
provides a map between the vertices of the triangulations, 
which are identified with $P^1(\mathbb Q)\subset S^1_{\infty}$. 
This map extends continuously to a homeomorphism of  $S^1_{\infty}$, 
which is piecewise-$PSL(2,\Z)$. 
This establishes an isomorphism between the Ptolemy group and 
the group of piecewise-$PSL(2,\Z)$ homeomorphisms of the circle. 

\vspace{0.2cm}
\noindent 
An explicit isomorphism with the group $T$ in the form introduced above 
was provided by Lochak and Schneps (see \cite{LS}). Send $\alpha$ to the flip 
$F_{a}$ of $(\tau_0, \vec{a_0})$ and $\beta$ to the 
element $((\tau_0, \vec{a_0}),(\tau_0, \vec{a_1}))$ of the Ptolemy group,
where $\vec{a_1}$ is the oriented edge in the basic triangle of  
the Farey triangulation $\tau_0$ next to $\vec{a_0}$.

\subsection{Quantum universal Teichm\"uller space}
Here and henceforth, for the sake of brevity, 
we will use the term tessellation instead of marked tessellation. 
For each tessellation $\tau$ let $E(\tau)$ be the set of its edges. 
We associate further a skew-symmetric 
matrix  $\varepsilon(\tau)$ with entries $\varepsilon_{ef}$,  for all $e,f\in E(\tau)$, 
as follows. If $e$ and $f$ do not belong to the same triangle of $\tau$ or $e=f$  then 
$\varepsilon_{ef}=0$. Otherwise, $e$ and $f$  are distinct edges belonging 
to the same triangle of $\tau$ and thus have a common vertex. 
We obtain $f$ by rotating $e$  in the plane along that vertex 
such that the moving edge  is sweeping out the  respective triangle of $\tau$. 
If we rotate  clockwisely then $\varepsilon_{ef}=+1$ and otherwise  
$\varepsilon_{ef}=-1$. 

\vspace{0.2cm}
\noindent 
The couple $(E(\tau), \varepsilon(\tau))$ is called a {\em seed} in \cite{FG}. 
Observe, that in this particular case seeds are completely determined by tessellations.

 \vspace{0.2cm}
\noindent 
Let $(\tau,\tau')$ be a flip $F_e$ in the edge $e\in E(\tau)$. Then the associated seeds 
$(E(\tau), \varepsilon(\tau))$  and  $(E(\tau'), \varepsilon(\tau'))$  
are obtained one from the other by a {\em mutation} in the direction $e$.  
Specifically, this means that   
there is an isomorphism $\mu_e:E(\tau)\to E(\tau')$ such that 
\[ \varepsilon(\tau')_{\mu_e(s)\mu_e(t)}=\left\{
\begin{array}{ll}
-\varepsilon_{st}, & \mbox{ \rm if } e=s \mbox{ \rm or } e=t, \\
\varepsilon_{st}, & \mbox{ \rm if } \varepsilon_{se}\varepsilon_{et} \leq 0, \\
\varepsilon_{st} + |\varepsilon_{se}|\varepsilon_{et},   
& \mbox{ \rm if } \varepsilon_{se}\varepsilon_{et} >0 
\end{array}
\right. 
\]
The map $\mu_e$ comes from the natural identification of the edges  of the two 
respective tessellations out of $e$ and $F_e(e)$. 

\vspace{0.2cm}
\noindent 
This algebraic setting appears in the description of the universal Teichm\"uller space 
${\mathcal T}$.  Its formal definition (see \cite{FG1,FG2})  is the  set of positive real 
points  of the cluster ${\mathcal X}$-space related to the set of seeds above. However, 
we can give a more intuitive description of it, following \cite{pe0}.  
Specifically,  ${\mathcal T}$ is the space of {\em all} marked tessellations  
(denoted ${\mathcal Tess}$ in \cite{pe0}). 
Each tessellation $\tau$ gives rise to a coordinate system 
$\beta_{\tau}:{\mathcal T}\to \R^{E(\tau)}$.  The real numbers 
$x_e=\beta_{\tau}(e)\in \R$ specify the amount of translation along the geodesic associated 
to the edge $e$  which is required when gluing together the two ideal triangles sharing that geodesic 
to obtain a given quadrilateral in the hyperbolic plane.  
These are called the shearing coordinates (introduced by Thurston and then 
considered by Bonahon, Fock et Penner) on the universal Teichm\"uller space and they provide a 
homeomorphism $\beta_{\tau}:{\mathcal T}\to \R^{E(\tau)}$. 
There is an explicit geometric formula (see also \cite{Fo,Fun}) for the shearing coordinates, as follows. 
Assume that the union of the two ideal triangles in ${\mathbb H}^2$ is the ideal quadrilateral 
of vertices $pp_0p_{-1}p_{\infty}$ and the common geodesic is $p_{\infty}p_0$. 
Then the respective shearing coordinate is the cross-ratio 
\[ x_e=[p,p_0,p_{-1},p_{\infty}]=\log\frac{(p_0-p)(p_{-1}-p_{\infty})}
{(p_{\infty}-p)(p_{-1}-p_0)}\]

\vspace{0.2cm}
\noindent 
Let $\tau'$ be obtained from $\tau$ by a flip $F_e$ and set $\{x'_f\}$ for the 
coordinates associated to $\tau'$.   The map  
$\beta_{\tau,\tau'}:R^{E(\tau')}\to \R^{E(\tau)}$ given  by 
\[ \beta_{\tau,\tau'}(x'_s)=\left\{\begin{array}{ll}
x_s-\varepsilon(\tau)_{se}\log(1+\exp(-{\rm sgn}(\varepsilon_{se})x_e)), & \mbox{\rm if } s\neq e,\\ 
-x_e, &  \mbox{\rm if } s= e\\ 
\end{array}
\right. 
\]
relates the two coordinate systems, namely $\beta_{\tau,\tau'}\circ \beta_{\tau'}=\beta_{\tau}$.  

\vspace{0.2cm}
\noindent 
These coordinate systems provide a contravariant functor $\beta: Pt\to {\rm Comm}$ 
from the Ptolemy groupoid $Pt$ to the category ${\rm Comm}$ of commutative topological 
$*$-algebras over $\C$. We associate to a tessellation $\tau$ the 
algebra $B(\tau)=C^{\infty}(\R^{E(\tau)}, \C)$ of smooth  complex valued 
functions on $\R^{E(\tau)}$, with the $*$-structure given by $*f=\overline{f}$. 
Further to any flip $(\tau,\tau')\in Pt$ one associates the map 
$\beta_{\tau,\tau'}:B(\tau')\to B(\tau)$. 

\vspace{0.2cm}
\noindent 
The matrices $\varepsilon(\tau)$ have a deep geometric meaning. 
In fact the  bi-vector  field 
\[ P_{\tau}=\sum_{e,f} \varepsilon(\tau)_{ef} \,\frac{\partial}{\partial x_e}\wedge 
\frac{\partial}{\partial x_f} \]
written here in the coordinates $\{x_e\}$ associated to $\tau$, 
defines a Poisson structure on ${\mathcal T}$ which is invariant by the action 
of the Ptolemy groupoid. The associated Poisson bracket is then given by the formula 
\[ \{x_e,x_f\}=\varepsilon(\tau)_{ef} \]

\vspace{0.2cm}
\noindent 
Kontsevich proved that  there is a canonical  formal quantization 
of a (finite dimensional) Poisson manifold. The universal Teichm\"uler space 
is not only a Poisson manifold but also endowed with a group action and 
our aim will be an equivariant quantization.  Chekhov, Fock  and Kashaev 
(see \cite{CF,K})  constructed an equivariant  quantization by means of 
explicit formulas. There are two ingredients in their approach. First, 
the Poisson bracket is given by constant coefficients, in any 
coordinate charts and  second,  the quantum (di)logarithm.

\vspace{0.2cm}
\noindent To any category $C$ whose morphisms are $\C$-vectors spaces one 
associates its projectivisation $PC$ having the same objects and new 
morphisms given by $Hom_{PC}(C_1,C_2)=Hom_C(C_1,C_2)/U(1)$, for any two objects 
$C_1,C_2$ of $C$. Here $U(1)\subset \C$ acts by scalar multiplication. 
A projective functor into $C$ is actually a functor into $PC$. 

\vspace{0.2cm}
\noindent 
Let now ${\rm A}^*$ be the category of topological $*$-algebras. 
Two functors $F_1, F_2:C\to {\rm A}^*$ essentially coincide if there exists a third functor $F$ and 
natural transformations $F_1\to F$, $F_2\to F$ providing dense inclusions 
$F_1(O)\hookrightarrow F(O)$ and $F_2(O)\hookrightarrow F(O)$, for any object $O$ of $C$.

\begin{definition}\label{quant}
A {\em quantization}  ${\mathcal T}^h$ of the universal Teichm\"uller space 
is a family of contravariant projective functors 
$\beta^h:Pt\to {\rm A}^*$ depending smoothly on the real parameter $h$ such that: 
\begin{enumerate}
\item The limit $\lim_{h\to 0} \beta^h= \beta^0$ exists and essentially coincide with the functor 
$\beta$. 
\item The limit $\lim_{h\to 0}[f_1,f_2]/h$ is defined and coincides with 
the Poisson bracket on ${\mathcal T}$. Alternatively, for each $\tau$ we have a 
$\C(h)$-linear (non-commutative) product structure $\star$ 
on the vector space $C^{\infty}(\R^{E(\tau)},\C(h))$ such that 
\[ f\star g= fg + h \{f,g\} + {o}(h)\]
where $\{f,g\}$ is the Poisson bracket on functions on ${\mathcal T}$ and 
$\C(h)$  denotes the algebra of smooth $\C$-valued functions on the real parameter $h$. 
\end{enumerate}
\end{definition}

\vspace{0.2cm}
\noindent 
We associate to each tessellation $\tau$ the {\em Heisenberg algebra} $H_{\tau}^h$ which is 
the topological $*$-algebra over $\C$ generated by the elements $x_e, e\in E(\tau)$ 
subjected to the relations
\[ [x_e,x_f]=2\pi i h \varepsilon(\tau)_{ef}, \;\,\; x_e^*=x_e \]
We define then $\beta^h(\tau)=H_{\tau}^h$. 

\vspace{0.2cm}
\noindent 
The quantization should associate a homomorphism 
$\beta^h((\tau,\tau')): H_{\tau'}^h\to H_{\tau}^h$ to each element 
$(\tau,\tau')\in Pt$. It actually suffices to consider the case when 
$(\tau,\tau')$ is the flip $F_e$ in the edge $e\in E(\tau)$. 
Let $\{x'_s\}$, $s\in E(\tau')$ be the generators of $H_{\tau'}^h$. 
We set then 
\[ \beta^h((\tau,\tau'))(x'_s)=\left\{\begin{array}{ll}
x_s -\varepsilon(\tau)_{se}\; \phi^h(-{\rm sgn}\left(\varepsilon(\tau)_{se})x_e\right), & \mbox{\rm if } 
s\neq e, \\
-x_s, & \mbox{\rm if } s= e \\
\end{array}
\right. 
\] 
Here $\phi^h$ is the {\em quantum logarithm} function, namely 
\[ \phi^h(z) = -\frac{\pi h}{2}\int_{\Omega}
\frac{\exp(-it z)}{{\rm sh}(\pi t) \, {\rm sh }(\pi h t)}  dt \]
where the contour $\Omega$ goes along the real axes from $-\infty$ 
to $\infty$ bypassing the origin from above. 

\vspace{0.2cm}
\noindent Some properties of the quantum logarithm are collected below: 
\[ \lim_{h\to 0}\phi^h(z)=\log\left(1+\exp(z)\right), \;\;  \phi^h(z)-\phi^h(-z)=z, \; \; \overline{\phi^h(z)}=\phi^h\left(\overline{z}\right), \;\; 
 \frac{\phi^h(z)}{h}=\phi^{1/h}\left(\frac{z}{h}\right)\]

\vspace{0.2cm}
\noindent 
A convenient way to represent this transformation graphically is to associate to a 
tessellation its dual binary tree embedded in $\H^2$ and to assign to each edge 
$e$ the respective generator $x_e$. Then the action of a flip reads as follows: 

\vspace{0.2cm}
\begin{center}
\includegraphics{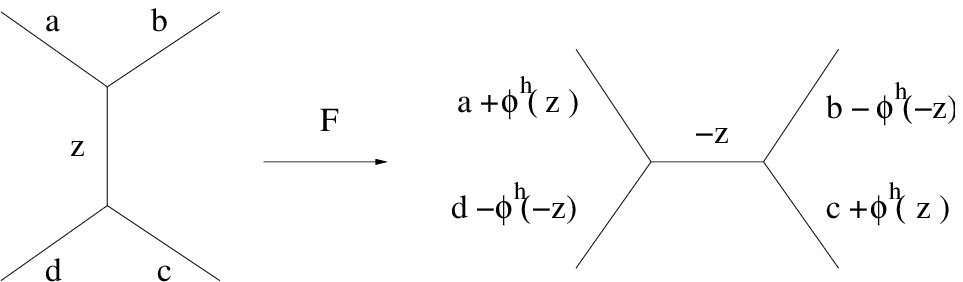}
\end{center}

\vspace{0.2cm}
\noindent We have then: 
\begin{proposition}[\cite{CF,FG3}]
The projective functor $\beta^h$ is well-defined and it is a quantization of the 
universal Teichm\"uller space $\mathcal T$. 
\end{proposition}

\vspace{0.2cm}
\noindent 
One proves that $\beta^h((\tau,\tau'))$ is independent on the decomposition of 
the element $(\tau,\tau')$ as a product of flips. 
In the  classical limit $h\to 0$ the quantum flip tends to the usual formula 
of the coordinates  change induced by a flip. 
Thus the first requirement in the definition \ref{quant} is fulfilled, and 
the second one  is obvious, from the defining relations in the Heisenberg algebra $H_{\tau}^h$.

\subsection{The dilogarithmic representation of $T$}
The subject of this section is to give a somewhat 
self-contained definition of the 
dilogarithmic representation of the group $T$. The case of general  cluster modular groupoids 
is developed in full detail in \cite{FG,FG3} and the group $T$ as a cluster modular groupoid 
is explained in \cite{FG2}.  

\vspace{0.2cm}
\noindent 
The quantization of a physical system in quantum mechanics should provide a Hilbert space 
and the usual procedure is to consider a Hilbert space representation of the algebra 
from definition \ref{quant}. This is formalized in the notion of 
representation of a quantum space. 

\vspace{0.2cm}
\noindent 
\begin{definition}\label{repre}
A projective $*$-representation of the quantized universal Teichm\"uller space ${\mathcal T}^h$,  
specified by the functor $\beta^h:Pt\to {\rm A}^*$,  
consists of the following data: 
\begin{enumerate}
\item A projective functor $Pt\to {\rm Hilb}$ to the category of Hilbert spaces. In particular, 
one associates a Hilbert space ${\mathcal L}_{\tau}$ to each tessellation $\tau$ and 
a unitary operator $K_{(\tau,\tau')}:{\mathcal L}_{\tau}\to {\mathcal L}_{\tau'}$, defined 
up to a scalar of absolute value 1.  
\item A $*$-representation $\rho_{\tau}$ of the Heisenberg algebra $H_{\tau}^h$ 
in the Hilbert space ${\mathcal L}_{\tau}$, such that the operators ${\mathbf K}_{(\tau,\tau')}$ 
intertwine the representations  $\rho_{\tau}$ and  $\rho_{\tau'}$ i.e.
\[   \rho_{\tau}(w)= {\mathbf K}_{(\tau,\tau')}^{-1} \rho_{\tau'}\left(\beta^h((\tau,\tau'))(w)\right)
{\mathbf K}_{(\tau,\tau')}, \,\,\, w\in H_{\tau}^h\]
\end{enumerate}
\end{definition}

\vspace{0.2cm}
\noindent 
{\em The classical Heisenberg $*$-algebra} $H$ is generated by $2n$ elements 
$x_s,y_s$, $1\leq s\leq n$ and relations 
\[ [x_s, y_s]=2\pi i \, h, \;\; [x_s,y_t]=0, \mbox{\rm if } s\neq t, \;\; [x_s, x_t]=[y_s,y_t]=0, \mbox{\rm for all } s, t\] 
with the obvious $*$-structure. 
The single irreducible integrable $*$-representation 
$\rho$ of $H$ makes it act   
 on the Hilbert space $L^2(\R^n)$ by means of the operators:  
\[ \rho(x_s)f(z_1,\ldots,z_n)= z_sf(z_1,\ldots,z_n) , \;\;\, \rho(y_s)=-2\pi i \, h \frac{\partial f}{\partial z_s} \]
{\em The Heisenberg algebras} $H_{\tau}^h$  are defined by commutation relations with 
constant coefficients and hence their representations can be constructed 
by selecting a Lagrangian subspace in the generators $x_s$   --  called a polarization  -- 
and let the generators act as linear combinations in the operators $\rho(x_s)$ and $\rho(y_s)$ 
above. 

\vspace{0.2cm}
\noindent 
The Stone von Neumann theorem holds true then for these algebras. 
Specifically, there exists an unique unitary irreducible Hilbert space representation 
of given central character that is integrable i.e. which can be integrated to the 
corresponding Lie group. Notice that  there exist in general also 
non-integrable unitary representations. 

\vspace{0.2cm}
\noindent 
In particular we obtain representations of $H_{\tau}^h$ and  $H_{\tau'}^h$. 
The uniqueness of the representation yields the existence of  an intertwinner ${\mathbf K}_{(\tau,\tau')}$
(defined up to a scalar) between the two representations. 
However,  nor the Hilbert spaces neither the representations $\rho_{\tau}$ are not canonical, 
as they depend on the choice of the polarization.

\vspace{0.2cm}
\noindent 
We will give below {\em the construction of a canonical representation} when
the quantized Teichm\"uller space  is replaced by its double. 
We need first to switch to another system of coordinates, 
coming from the cluster  ${\mathcal A}$-varieties. 
Define, after Penner (see \cite{pe0}) 
the  {\em universal decorated Teichm\"uller space} ${\mathcal A}$ to be the space 
of all marked tessellations endowed with one horocycle for each vertex (decoration).   
Alternatively (see \cite{FG1}), 
${\mathcal A}$ is the set of positive real points of the 
cluster ${\mathcal A}$-space related to the previous set of seeds. 

\vspace{0.2cm}
\noindent 
Each tessellation $\tau$  yields a coordinate system 
$\alpha_{\tau}:{\mathcal A}\to \R^{E(\tau)}$ which associates to the edge 
$e$ of $\tau$ the coordinate $a_e=\alpha_{\tau}(e)\in \R$. 
The number $\alpha_{\tau}(e)$  is the algebraic distance between 
the two horocycles on $\H^2$ centered at vertices of $e$, 
measured along the geodesic associated to $e$.  These are the so-called {\em lambda} 
coordinates of Penner. 

\vspace{0.2cm}
\noindent 
There is a canonical  map $p:{\mathcal A}\to {\mathcal T}$  (see \cite{pe0}, Proposition 3.7 and 
\cite{FG1}) 
such that, in the coordinate systems induced by a tessellation $\tau$, the corresponding map 
$p_{\tau} :\R^{(E(\tau)}\to \R^{E(\tau)}$ is given by 
\[ p_{\tau}\left(\sum_{t\in E(\tau)} \varepsilon(\tau)_{st}a_t\right)=x_s\]

\vspace{0.2cm}
\noindent 
Let  $(\tau,\tau')$ be the flip on the edge $e$ and set $a'_s$ be the coordinates system 
associated to $\tau'$. Then the flip induces the  following change of coordinates:  
\[ \alpha_{\tau,\tau'}(a_s)=\left\{\begin{array}{ll}
a_s, & \mbox{\rm if } s\neq e \\
-a_e + \log\left(\exp\left(\sum_{t;\varepsilon(\tau)_{et}>0}
\varepsilon(\tau)_{et} a_t\right)+ \exp\left(-\sum_{t;\varepsilon(\tau)_{et}<0}
\varepsilon(\tau)_{et} a_t\right)\right), & \mbox{\rm if } s= e \\
\end{array}
\right.
\]
It can be verified that $p_{\tau}$ are compatible with the action of the 
Ptolemy groupoid on the respective coordinates.

\vspace{0.2cm}
\noindent 
{\em The vector space} ${\mathcal L}_{\tau}$ is defined as the space of 
square integrable functions  with finite dimensional support 
on ${\mathcal A}$ with respect to the $\alpha_{\tau}$ coordinates 
i.e. the functions $f:\R^{E(\tau)}\to \C$,  with support contained into 
some $R^F\times \{0\}\subset R^{E(\tau)}$, for some finite subset 
$F\subset E(\tau)$. The coordinates on $\R^{E(\tau)}$ are the $a_e, e\in E(\tau)$. 
The function $f$ is square integrable if 
\[ \int_{\R^F} |f|^2 \bigwedge _{e\in F}da_e  < \infty \]
for any such $F$ as above. Let $f, g \in {\mathcal L}_{\tau}$. Then let 
 $\R^F\times \{0\}$ contain the intersection of their supports. 
Choose $F$ minimal with this property. Then the scalar product 
\[ \langle f, g \rangle= \int_{\R^F} f(a)\overline{g(a)} \bigwedge _{e\in F}da_e \]
makes ${\mathcal L}_{\tau}$ a Hilbert space.

\vspace{0.2cm}
\noindent 
To define the intertwining operator ${\mathbf K}$ we set now: 
\[ G_e((a_s)_{s\in F})=
\int \exp\left(\int_{\Omega}\frac{\exp(it \sum_{s\in F}\varepsilon(\tau)_{es}a_s) \sin(tc)}
{2i {\rm sh}(\pi t){\rm sh}(\pi h t)} \frac{dt}{t} + \frac{c}{\pi i h}\left(\sum_{s; \varepsilon(\tau)_{es}<0}\varepsilon(\tau)_{es}a_s 
+a_e\right)\right) dc \]
The key ingredient in the construction of this function is the 
{\em quantum dilogarithm} (going back to  Barnes (\cite{Bar}) and 
used by Baxter (\cite{Bax}) and Faddeev (\cite{Fad})): 
\[ \Phi^h(z) = \exp\left(-\frac{1}{4}\int_{\Omega}
\frac{\exp(-it z)}{{\rm sh}(\pi t) \, {\rm sh }(\pi h t)}  \frac{dt}{t} \right)\]
where the contour $\Omega$ goes along the real axes from $-\infty$ 
to $\infty$ bypassing the origin from above. 

\vspace{0.2cm}
\noindent Some properties of the quantum dilogarithm are collected below: 
\[ 2\pi i h d \log \Phi^h(z)=\phi^h(z), \; \;  \lim_{\Re z\to -\infty}\Phi^h(z)=1\]
\[ \lim_{h\to 0}\Phi^h(z)/\exp(-{\rm Li}_2(-\exp(z)))=2\pi i h, \, \, \mbox{\rm where } 
{\rm Li}_2(z)=\int_0^z\log(1-t)dt \]
\[ \Phi^h(z)\Phi^h(-z)=\exp\left(\frac{z^2}{4\pi i h}\right) \exp\left(-\frac{\pi i}{12}(h+h^{-1})\right), \;\; 
 \overline{\Phi^h(z))}=\left(\Phi^h\left(\overline{z}\right)\right)^{-1}, \;\; \Phi^h(z)=\Phi^{1/h}\left(\frac{z}{h}\right)\]

\vspace{0.2cm}
\noindent 
Let now $f\in {\mathcal L}_{\tau}$, namely some  $f:\R^F\times \{0\}\to \C$. 
Consider $(\tau,\tau')$ be the flip $F_e$ on the edge $e$. 
Let $a_s, s\in F$ be the coordinates in $\R^F$. If $e\not\in F$ then we set 
\[ {\bf K}_{(\tau,\tau')}=1\]
 If $e\in F$ then the coordinates associated to $\tau'$ are $a_s, s\neq e$ and $a'_e$. 
Set then 
\[ ({\bf K}_{(\tau,\tau')} f)((a_s, _{s\in F, f\neq e}, a'_e)= 
\int G_e((a_s)_{s\in F, s\neq e}, a_e+a'_e) f((a_s)_{s\in F}) da_s \]

\vspace{0.2cm}
\noindent 
The last piece of data is {\em the representation of the Heisenberg algebra} 
$H_{\tau}^h$ in the Hilbert space ${\mathcal L}_{\tau}$. We can 
actually  do better, namely to enhance the space with a bimodule structure. 
Set 
\[ \rho^-_{\tau} (x_s) =-\pi i h \frac{\partial}{\partial a_s} + \sum_{t}\varepsilon(\tau)_{st} a_t \]
\[ \rho^+_{\tau} (x_s)=\pi i h \frac{\partial}{\partial a_s} + \sum_{t}\varepsilon(\tau)_{st} a_t \]
Then $\rho^-_{\tau} $ gives a left module and $\rho^+_{\tau} $ a  right module structure on 
${\mathcal L}_{\tau}$ and the two actions commute. 
We have then: 

\begin{proposition}[\cite{CF,FG3,FG}]
The data $({\mathcal L}_{\tau}, \rho^{\pm}_{\tau}, {\mathbf K}_{(\tau,\tau')})$ 
is a projective $*$-representation of the  quantized universal Teichm\"uller space.   
\end{proposition}

\vspace{0.2cm}
\noindent
We call it the dilogarithmic representation of the Ptolemy groupoid. 
The proof of this result is given in \cite{FG} and  a particular case is explained 
with lots of details  in \cite{Go}. 

\vspace{0.2cm}
\noindent
The last step in our construction is to observe that a 
representation of the Ptolemy groupoid 
$Pt$ {\em induces a representation of the Ptolemy-Thompson group} 
$T$ by means of 
an identification of the Hilbert spaces ${\mathcal L}_{\tau}$ for all $\tau$.

\vspace{0.2cm}
\noindent
Projective representations are equivalent to representations of 
central extensions by  means of the following well-known procedure.
To a general group $G$,  Hilbert space $V$ and homomorphism 
$A:G\to PGL(V)$ we can associate a central extension $\widetilde{G}$ 
of $G$ by $\C^*$ which resolves the projective representation $A$ 
to a linear representation $\widetilde{A}:\widetilde{G}\to GL(V)$. 
The extension $\widetilde{G}$ is the pull-back on $G$ of the 
canonical central $\C^*$-extension $GL(V)\to PGL(V)$. 

\vspace{0.2cm}
\noindent
However the central extension which we consider here is a subgroup 
of the $\C^*$-extension defined above, obtained by using 
a particular section over $G$. Let us 
write $G=F/R$ as the quotient of the free group $F$ by the normal 
subgroup $R$ generated by the relations.  Then our data 
consists in a homomorphism $\overline{A}:F\to GL(V)$ with 
the property that $\overline{A}(r)\in\C^*$, for each 
relation $r\in R$, so that $\overline{A}$ induces $A:G\to PGL(V)$. 
This data will  be called 
{\em an almost-linear representation}, in order to distinguish 
it from a  projective representation of $G$. 

\vspace{0.2cm}
\noindent
The {\em central extension $\widehat{G}$ of $G$ associated to 
$\overline{A}$} 
is $\widehat{G}=F/(\ker \overline{A}\cap R)$, namely the smallest  central 
extension  of $S$ resolving the projective representation $A$ to a linear 
representation compatible with $\overline{A}$. Then 
$\widehat{G}$ is a central extension of $G$ by the subgroup 
$\overline{A}(R)\subset \C^*$ and hence it is naturally a subgroup 
of $\widetilde{G}$. In other terms $\overline{A}$ determines a 
projective representation $A$ and a section over $G$ whose 
associated 2-cocycle takes values in $\overline{A}(R)$ and which 
describes the central extension $\widehat{G}$.

\vspace{0.2cm}
\noindent
Now, the intertwinner functor $\mathbf K$  is actually an almost-linear  
representation (in the obvious sense) of the 
Ptolemy groupoid and thus induces an  
almost-linear representation of the Ptolemy-Thompson group $T$ 
into the unitary group. 
We can extract from \cite{FG} the following results  (see also the equivalent 
construction at the level of Heisenberg  algebras in \cite{BBL}):

\begin{proposition}\label{dilogcar}
The dilogarithmic almost-linear 
representation $\mathbf K$ has  the following properties: 
\begin{enumerate}
\item Images of disjoint flips in $\widehat{T}$  commute with each other; 
\item The square of a flip is the identity;
\item The composition of the  lifts of the five flips from the 
pentagon relation below is $\exp(2\pi i h)$ times the symmetry permuting the 
two edges coordinates. 

\begin{center}
\includegraphics{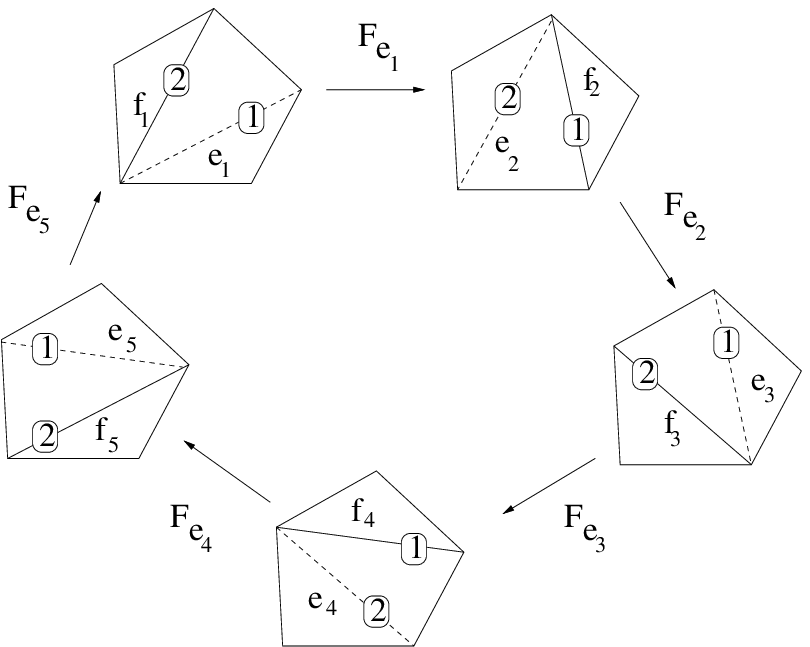}
\end{center}

\end{enumerate}
\end{proposition}
\begin{proof}
The first condition is that images by $K$ of flips on disjoint edges 
should commute. This is obvious by the explicit formula for $K$. 
The second and third conditions are proved in \cite{FG}. 
\end{proof}

\vspace{0.1cm}
\noindent 
Therefore the image by $\mathbf K$ of relations of the Ptolemy groupoid 
into $\C^*$ is the subgroup $U$ generated by $\exp(2\pi i h)$.  
We can view  the pentagon relation in the Ptolemy-Thompson 
group $T$ as a pentagon relation in the Ptolemy groupoid $Pt$. Thus  
the image by $\mathbf K$ of relations of the Ptolemy-Thompson group $T$ 
into $\C^*$  is  also the subgroup $U$. In particular the associated 
2-cocycle takes values in $U$.  
If $h$ is a formal parameter or an irrational 
real number we obtain then a 2-cocycle with values in $\Z$. 

\begin{definition}
The dilogarithmic central extension 
$\widehat{T}$ is the central extension of $T$ by $\Z$  
associated to the dilogarithmic almost-linear representation $\mathbf K$ 
of $T$, or equivalently, to the previous 2-cocycle. 
\end{definition}

\subsection{Identifying the two central extensions of $T$}

\vspace{0.2cm}
\noindent 
The main result of this section is the following: 

\begin{proposition}
The dilogarithmic extension $\widehat{T}$ is identified to $T^*_{\rm ab}$. 
\end{proposition}
\begin{proof}
The main step is to translate the properties of the dilogarithmic 
representation of the Ptolemy groupoid in terms of the Ptolemy-Thompson group.  
Since $\widehat{T}$ is a central extension 
of $T$ it is generated by the lifts $\al,\be$ of $\alpha$ and $\beta$ 
together with the generator $z$ of the center. 
Let us see what are the relations arising in the 
group $\widehat{T}$.  According to proposition \ref{dilogcar} 
lifts of disjoint flips should commute. By a simple computation we can 
show that 
the elements $\beta\alpha\beta$, $\alpha^2\beta\alpha\beta\alpha^2$ and 
$\alpha^2\beta\alpha^2\beta\alpha\beta\alpha^2\beta^2\alpha^2$ act as 
disjoint flips on the Farey triangulation. 
In particular we have the relations 
\[  [\be\al\be,\al^2\be\al\be\al^2]=
 [\be\al\be,\al^2\be\al^2\be\al\be\al^2\be^2\al^2]=1\]
satisfied in $\widehat{T}$. Moreover, by construction we have also 
\[\be^3=\al^4=1\]
meaning that the $\al$ is still periodic of order $4$ while $\be$ is not deformed.

\vspace{0.2cm}
\noindent 
Eventually the only non-trivial lift of relations comes from the 
pentagon relation $(\be\al)^5$.  The element $(\be\al)^5$ is actually 
the permutation of the two edges in the pentagon times the composition of the five flips.  
The pentagon equation is not anymore satisfied 
but proposition \ref{dilogcar}  shows that the dilogarithmic image of  
$(\be\al)^5$ is a scalar operator. 
Since $z$ is the generator of the 
kernel $\Z$ of $\widehat{T}\to T$  
it follows that the the lift of the pentagon equation from $T$ to 
$\widehat{T}$ is given by 
\[ (\be\al)^5=z\]
According to proposition \ref{abelcar} all relations presenting 
$T^*_{\rm ab}$ are satisfied  in $\widehat{T}$. Since $\widehat{T}$ is a 
nontrivial  central extension of $T$ by $\Z$ 
it follows that the groups are isomorphic.  
\end{proof}

\begin{remark}
The key-point in the proof above is that all pentagon relations in $Pt$ 
are transformed in a single pentagon relation in $T$ and thus 
the scalars associated to the pentagons in $Pt$  should be the same.   
\end{remark}

\section{Classification of central extensions of the group $T$}
\subsection{The family  $T_{n,p,q,r}$ of central extensions}
Our main concern here is to  identify the cohomology classes  
of all central extensions of $T$ in $H^2(T)$. 
Before doing that we consider a series of 
central extensions $T_{n,p,q,r,s}$ of $T$ by $\Z$, having properties 
similar to those of  $\widehat{T}$. 

\begin{definition}
The group $T_{n,p,q,r,s}$, is presented by  the generators  
$\als, \bls, z$ and the relations: 
\[ (\bls\als)^5=z^n \]  
\[ \als^4=z^p \]
\[ \bls^3=z^q\]
\[ [\bls\als\bls, \als^2\bls\als\bls\als^2]=z^r\]
\[[\bls\als\bls, \als^2\bls\als^2\bls\als\bls\als^2\bls^2\als^2]=z^s\]
\[ [\als,z]=[\bls,z]=1\]
Let us denote  $T_{n,p,q,r}=T_{n,p,q,r,0}$ and 
$T_{n,p,q}=T_{n,p,q,0,0}$. 
\end{definition}

\vspace{0.2cm}
\noindent 
According to \cite{FK2} we can identify $\widehat{T}$ with 
$T_{1,0,0}$. In fact the group $T^*$ is split over 
 the smaller Thompson group $F\subset T$ and thus $\widehat{T}$ is 
split over $F$. Further $F$ is generated by the elements 
$\beta^2\alpha$ and $\beta\alpha^2$ and thus relations of $F$ are precisely 
given by the commutation relations above. 
Thus the last two relations hold true, while 
$z$ is central and thus $\widehat{T}$ is given by the presentation above. 

\begin{remark} 
We considered in \cite{FK2} the twin group $T^{\sharp}$ 
and gave a presentation 
of it.  Then, using a similar procedure there is a group obtained from 
$T^{\sharp}$ 
by abelianizing the kernel $B_{\infty}$, which is identified  
actually  to $T_{3,1,0}$. 
\end{remark}

\begin{proposition}\label{secondcom}
Central extensions of $T$ by $\Z$ are exhausted by the 
set of extensions $T_{n,p,q,r}$. 
\end{proposition}
\begin{proof}
Consider the relations in $T_{n,p,q,r,s}$ other than  
\[[\bls\als\bls, \als^2\bls\als^2\bls\als\bls\als^2\bls^2\als^2]=z^s\]
It suffices to see that these relations already force 
\[[\bls\als\bls, \als^2\bls\als^2\bls\als\bls\als^2\bls^2\als^2]=1\]
The commutator 
$[\beta\alpha\beta, \alpha^2\beta\alpha^2\beta\alpha\beta
\alpha^2\beta^2\alpha^2]$ is the trivial element of $T$ and thus 
it leads by means of Hopf theorem to a 2-cycle on $T$, given by 
the formula:
\[ \delta=(\beta\alpha\beta, \alpha^2\beta\alpha^2\beta\alpha\beta
\alpha^2\beta^2\alpha^2)- (\alpha^2\beta\alpha^2\beta\alpha\beta
\alpha^2\beta^2\alpha^2,\beta\alpha\beta)\]   
However the commutator above can be written as a commutator in the 
subgroup $F$ as we already remarked that 
\[ [\beta\alpha\beta, \alpha^2\beta\alpha^2\beta\alpha\beta
\alpha^2\beta^2\alpha^2]=[AB^{-1}, A^{-2}BA^2]\]
where $A=\beta\alpha^2$ and $B=\beta^2\alpha$ are the generators of $F$. 
Moreover the last commutator defines the 2-cycle 
\[ \epsilon=(AB^{-1}, A^{-2}BA^2)-(A^{-2}BA^2,AB^{-1})\]
in $H_2(F)$ and the inclusion $i:F\subset T$ sends the class  of 
$[\epsilon]$ into $i_*([\epsilon])= [\delta]$. 
However,  it is known  from \cite{GS} that 
$[\epsilon]$ is in the kernel of $i_*:H_2(F)\to H_2(T)$ and thus 
$[\delta]=0$. This shows that 
$\langle c_{T_{n,p,q,r,s}},[\epsilon]\rangle=0$, 
where 
$\langle, \rangle:H^2(T)\times H_2(T)\to \Z$ is the obvious pairing. 
Thus $T_{n,p,q,r,s}$ splits over the subgroup generated by 
$\beta\alpha\beta$ and $\alpha^2\beta\alpha^2\beta\alpha\beta
\alpha^2\beta^2\alpha^2$.
Thus up to changing each one of 
the lifts $\als,\bls\in T_{n,p,q,r,s}$ of $\alpha,\beta$ by 
a central factor we have 
\[[\bls\als\bls, \als^2\bls\als^2\bls\als\bls\als^2\bls^2\als^2]=1\]
However, if we chose another lifts then central factors will cancel each other in 
the commutator above and thus the identity above holds  true for any choice of the lifts. In particular 
we have $s=0$.   
\end{proof}

\vspace{0.2cm}
\noindent 
The aim of this chapter is to prove theorem \ref{exte0}, namely to compute  
the class $c_{T_{n,p,q,r}}\in H^2(T)$ of the extension $T_{n,p,q,r}$. 
  We denote by $\chi(n,p,q,r)$ the coefficient of $\chi$ and 
$\alpha(n,p,q,r)$ the coefficient of $\alpha$ in $c_{T_{n,p,q,r}}$.

\subsection{Computing $\alpha(n,p,q,r)$}
\begin{proposition}\label{alpha}
We have $\alpha(n,p,q,r)=r$.
\end{proposition}
\begin{proof}
Since the commutator $[\beta\alpha\beta,\alpha^2\beta\alpha\beta\alpha^2]$ is the  
identity in $T$  it gives rise to a 2-cycle in homology given by 
\[ \mu=(\beta\alpha\beta, \alpha^2\beta\alpha\beta\alpha^2)-
(\alpha^2\beta\alpha\beta\alpha^2,\beta\alpha\beta)\]
and hence representing a class $[\mu]\in H_2(T)$. 
As in the proof of Proposition \ref{secondcom} $[\mu]=i_*([\eta])$ where 
\[ \eta=(AB^{-1}, A^{-1}BA)-(A^{-1}BA,AB^{-1})\]
is the 2-cycle on $F$ associated to the commutator 
$[AB^{-1}, A^{-1}BA]$.

\vspace{0.2cm}\noindent 
We claim that  
\begin{lemma}
$\langle \alpha, [\mu]\rangle=1$, where 
$\langle, \rangle:H^2(T)\times H_2(T)\to \Z$ is the obvious pairing. 
\end{lemma}
\begin{proof}
Consider the discrete Godbillon-Vey 2-cocycle 
$\overline{gv}:T\times T\to \Z$ defined by the formula 
\[ \overline{gv}(g,h)=\sum_{x\in S^1}\det \left (
\begin{array}{cc}
\log_2h'_r(x) & \log_2(g\circ h)'_r(x) \\
h''(x)       & (g\circ h)''(x)
\end{array}
\right)
\]
Here $\gamma''$ stands for 
\[ \gamma''(x)=\log_2 \gamma'_r(x)- \log_2 \gamma'_l(x)\] 
where $\gamma'_r, \gamma'_l$ are the right respectively left derivatives 
of $\gamma$. Notice that the derivative of $\gamma\in T$ is a locally 
constant function having only finitely many discontinuity points.

\vspace{0.2cm}\noindent 
It is well-known (see \cite{GS}) that the 2-cocycle $\overline{gv}$ 
represents $2\alpha$ in cohomology. 
A direct computation using this cocycle shows that 
$\overline{gv}(\mu)=2$ and hence proving the claim. 
\end{proof}

\vspace{0.2cm}\noindent 
We have
\[ \langle c_{T_{n,p,q,r}}, i_*([\eta])\rangle=
\langle c_{T_{n,p,q,r}}|_{F}, [\eta]\rangle=
\langle c_{F_{r}}, [\eta]\rangle \]
where $F_r$ is the central extension of $F$ given by 
\[ F_r=\langle A,B,z; [AB^{-1}, A^{-1}BA]=z^r, \, 
[AB^{-1}, A^{-2}BA^2]=1, [z,A]=[z,B]=1\rangle \]
According to the arguments from the proof of 
Proposition \ref{secondcom} these extensions  exhaust the 
set of all central extensions of $F$ that arise from 
central extensions of $T$. 

\vspace{0.2cm}\noindent 
On the other hand we have 
\begin{eqnarray*} \langle c_{T_{n,p,q,r}}, [\mu])\rangle &=&
\langle \alpha(n,p,q,r)\alpha +\chi(n,p,q,r)\chi, i_*[\eta]\rangle=
\langle \alpha(n,p,q,r)\alpha|_{F} +\chi(n,p,q,r)\chi|F, [\eta]\rangle\\
&=&
\alpha(n,p,q,r)\langle \alpha|F,[\eta]\rangle=
\alpha(n,p,q,r)\end{eqnarray*}
because the image of $\chi$ in $H^2(F)$ vanishes.

\vspace{0.2cm}\noindent 
Observe that $F_r$ does not depend on $n,p,q$ and thus its class does not depend either. This shows that $\alpha(n,p,q,r)$ is a function on $r$. 
This function is computed through the choice of a setwise section 
in the projection $F_r\to F$, which is further evaluated to a fixed  
2-cycle $\eta$ and thus it should be a linear function  with integer 
coefficients and without constant 
terms. This implies that $\alpha(n,p,q,r)=\lambda r$, for $\lambda\in\Z$.   

\vspace{0.2cm}\noindent 
Eventually by choosing explicit lifts of the elements $AB^{-1}$ and $A^{-1}BA$ 
(and their products in whatever order) in $F_r$ we can evaluate 
$\langle c_{F_{r}}, [\eta]\rangle=r$. 

\end{proof}

\begin{corollary}
The extensions $T_{n,p,q,r}$ describe all central extensions of $T$
by $\Z$.

\vspace{0.2cm}\noindent 
The extensions $T_{n,p,q,0}$ describe all central extensions of $T$
by $\Z$  associated to multiples of the Euler class, and all 
of them split over $F$. 
\end{corollary}
\begin{proof}
The groups $T_{n,p,q,0}$ are extensions of $T$ by a cyclic group. 
Moreover any central extension of $T$ has the form $T_{n,p,q,r}$. 
The formula for $\alpha(n,p,q,r)$ shows that $r=0$ if we want that 
its class be a multiple of the Euler class. 
Further the set of central extensions of $T$ by $\Z$ 
up to isomorphism over $T$ are in one-to-one correspondence 
with the classes of $H^2(T)$. In particular any multiple of $\chi$ 
should be realized by one extension by $\Z$. Now, for given 
$m\in\Z$ there exists only one  isomorphism type 
of $T_{n,p,q,0}$ for which $n,p,q$ are solutions of 
$12n-15p-20q=m$. This follows from the uniqueness of solutions up to 
equivalence $p\sim p+4x$, $q\sim q+3y$, $n\sim n+5x+5y$.  
Thus the respective extension must have kernel $\Z$.  
\end{proof}

\subsection{The Greenberg-Sergiescu extension ${\mathcal A}_T$}
We know the value of $\chi(n,p,q,r)$ for $r=0$, but in order to find the 
coefficient of $r$ we need to know explicitly central extensions 
with nontrivial Godbillon-Vey class. Fortunately, Greenberg and Sergiescu 
in \cite{GrS} had constructed such an extension of $T$ by $B_{\infty}$. 
The main difficulty in analyzing this extension comes from the fact that there are 
several perspectives for analyzing the group $T$, 
either as a group of dyadic 
piecewise affine homeomorphisms 
of $S^1$ or else as a group $PPSL(2,\Z)$ of piecewise $PSL(2,\Z)$ automorphisms of the circle at infinity. 
If we plug in the discrete Godbillon-Vey then the formulas 
from \cite{GS,GrS} use the first point of view. If we are seeking for 
the mapping class group perspective it is the second point of view which 
is manifest and there is no direct relationship between this and 
the former one. The key-point in the calculations below is to pass from one 
perspective to the other. We have therefore to give a detailed account 
of the group ${\mathcal A}_T$, following Kapoudjian and Sergiescu \cite{KS}. 

\vspace{0.2cm}\noindent 
We will use the mapping class group description of the Ptolemy-Thompson 
group but we will enlarge the surface. We follow closely 
\cite{KS,FK2}.


\vspace{0.2cm}\noindent 
The surface $D$ occurs in the process of understanding the almost action 
of $T$ on the infinite binary tree. Recall that an almost automorphism 
of a tree is a map sending the complement of a finite tree 
isomorphically on the complement of a finite tree. The action of $T$ 
acting as mapping class group of $D$ induces an almost action on the 
binary tree. This point of view emphasize the  realization of 
$T$ as the group $PPSL(2,\Z)$.  

\vspace{0.2cm}\noindent 
 We can  think of $T$ as the group of dyadic piecewise affine homeomorphisms 
of $S^1$. This will lead to another, more subtle, 
way to construct an action of $T$ on a regular tree. 
Let ${\mathcal T}$ be the rooted binary tree with one 
finite leaf, obtained by splitting one edge of the usual binary tree 
at some vertex (the root) and attaching one more edge. 
Label the vertices and edges (see \cite{CFP})
of ${\mathcal T}$  inductively as follows. The leaf is labeled $0\sim 1$ and 
the edge joining it with the root by $[0,1]$. The root is 
labeled $\frac{1}{2}$,  its left descending edge 
is labeled $[0,\frac{1}{2}]$ and the right descending one $[\frac{1}{2},1]$. 
Further if an edge is labeled $[\frac{k}{2^n},\frac{k+1}{2^n}]$, then 
its bottom vertex is labeled $\frac{2k+1}{2^{n+1}}$ and the two edges 
issued from it are labeled, the 
left one $[\frac{k}{2^n},\frac{2k+1}{2^{n+1}}]$ 
and the right one $[\frac{2k+1}{2^{n+1}},\frac{k+1}{2^n}]$ respectively.

\begin{center}
\includegraphics{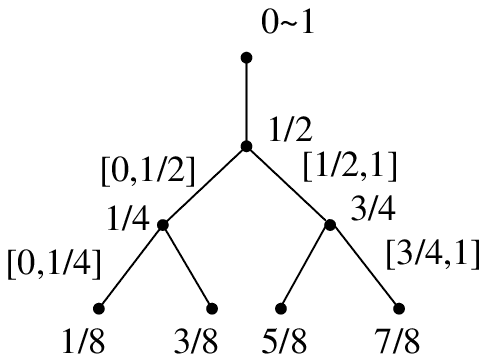}
\end{center}
  
\vspace{0.2cm}\noindent 
The almost action of $T$ on ${\mathcal T}$ could be read off 
from the description of elements of $T$ as pairs of stable binary trees. 
Alternatively, we identify elements of $T$ as dyadic piecewise 
affine homeomorphisms of $S^1=[0,1]/0\sim 1$. The action of such a 
homeomorphism induces a bijection of the set of vertices of ${\mathcal T}$ 
(identified to their labels). This bijection is an almost automorphism 
of the tree ${\mathcal T}$.

\vspace{0.2cm}\noindent 
Consider next the tree ${\mathcal ET}$ obtained from ${\mathcal T}$   
by adjoining a pending line (with infinitely many vertices on it) to each 
vertex of ${\mathcal T}$.

\begin{center}
\includegraphics{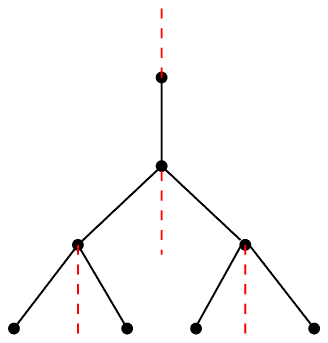}
\end{center}
   
\vspace{0.2cm}\noindent 
There is an obvious extension of the almost action of $T$ from 
${\mathcal T}$ to ${\mathcal ET}$. However there exists a more 
interesting one as it was discovered in \cite{GrS}.

\begin{definition}
Let $\Q\subset S^1$ denote the set of images of dyadic numbers. 
Assume that $T$ is identified with the subgroup of ${\rm Homeo}^+(S^1)$ 
of piecewise linear homeomorphisms. 
A {\em cocycle} is a map $K:T\times \Q\to \Z$  satisfying the following conditions: 
\begin{enumerate}
\item  For any $\gamma\in T$,  $K(\gamma, x)$ vanishes for 
all but finitely many points $x\in \Q$ and 
\[ \sum_{x\in \Q}K(\gamma,x)=0\]
\item $K(\gamma\delta, w)=K(\gamma, \delta(w))+K(\delta,w)$, for all $\gamma,\delta\in T$, $w\in \Q$. 
\end{enumerate}
\end{definition}

\vspace{0.2cm}\noindent 
For any cocycle $K:T\times \Q\to \Z$ we can associate such an action. 
Let $\gamma\in T$. Then $\gamma$ induces a bijection denoted by the 
same letter between the vertices of the rooted tree ${\mathcal T}$, which 
were identified with $\Q$. 
This bijection induces an almost automorphism of the tree ${\mathcal T}$. 
Moreover, let $n=\max_{x\in \Q} |K(\gamma,x)|$. 
If $v\in \Q$ is a vertex of $\mathcal T$ let $f_v$ denote the 
pending line at $v$ and $f^{\geq n}_v\subset f_v$ be the subtree of 
those points at distance at least $n$ from $v$. 
We define the almost action of $\gamma$ on $\mathcal ET$ as being 
the unique isometric bijection from 
$f^{\geq n}_{v}$ to $f^{\geq n}_{\gamma(v)}$.

\begin{definition}
The enhanced ribbon tree $ED$  is obtained 
by thickening $\mathcal ET$ in the plane and $ED^*$ by  puncturing 
$ED$ along the vertices of the pending lines. 
Notice that the vertices of $\mathcal T$ are not among the 
punctures.  

\begin{center}
\includegraphics{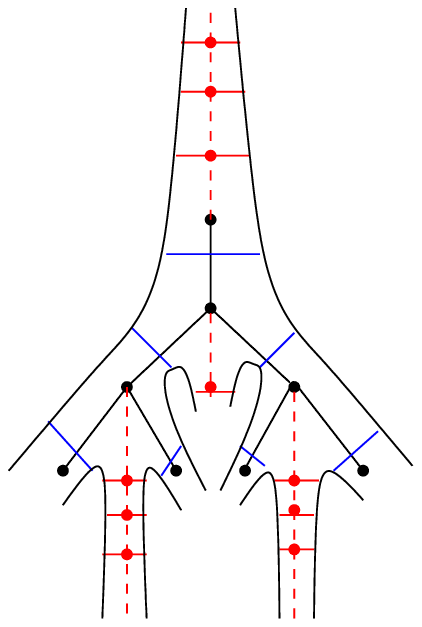}
\end{center}

\vspace{0.2cm}\noindent 
The canonical rigid structures of $ED$ and $ED^*$ are those 
from the picture above which 
decomposes $ED$ into  (punctured) octagons centered at vertices of 
$\mathcal T$ and  (punctured) squares along the pending lines at 
vertices of $\mathcal T$.

\vspace{0.2cm}\noindent Admissible sub-surfaces of $ED$ or $ED^*$ are
connected finite unions of elementary pieces.  
\end{definition}

\vspace{0.2cm}\noindent 
Given a cocycle $K$ one finds  that the almost action 
of $T$ on $\mathcal ET$ induces an embedding of the group $T$ into 
the asymptotic mapping class group 
${\mathcal M}(ED)$ with the given structure. Specifically, define 
the content of an admissible subsurface $\Sigma\subset ED$ to be 
the number of squares it contains. Then $T$ is the group of 
mapping classes that preserve the content i.e. those mapping classes 
of homeomorphisms $\varphi$ for which 
$\varphi(\Sigma)$ and $\Sigma$ have the same content for any admissible 
$\Sigma$. 

\vspace{0.2cm}\noindent 
We have an obvious exact sequence 
\[ 1\to B_{\infty}\to {\mathcal M}(ED^*)\to {\mathcal M}(ED)\to 1 \]
Using the embedding $\iota_K:T\hookrightarrow {\mathcal M}(ED)$ 
we can restrict  ${\mathcal M}(ED^*)$ at $\iota_K(T)$. This restriction 
is the group ${\mathcal A}_{T,K}$, which fits also in the 
exact sequence 
\[ 1\to B_{\infty}\to {\mathcal A}_{T,K}\to T\to 1 \]

\vspace{0.2cm}\noindent 
The main example of a nontrivial cocycle is the one associated to the 
discrete Godbillon-Vey class. Specifically, identify $T$ with the 
group of dyadic piecewise affine homeomorphisms of $[0,1]/0\sim 1$. 
For any $\gamma\in T$ and $v\in \Q$ we set 
\[ K(\gamma,v)=\gamma''(v)=\log_2 \gamma'_r(v)- \log_2 \gamma'_l(v)\] 
where $\gamma'_r, \gamma'_l$ are the right respectively left derivatives 
of $\gamma$.  It is well-known (\cite{GS,GrS})
that $K(\gamma,v)$ is a cocycle. 

\vspace{0.2cm}\noindent 
The extension ${\mathcal A}_{T,K}$ obtained when $K$ is the 
Godbillon-Vey cocycle above is simply denoted ${\mathcal A}_T$.

\begin{remark}
The definition from \cite{KS} was slightly different because it used $n+1$ 
instead of $n$ in the definition  and we have punctures at 
the vertices of $\mathcal T$; in particular the vertex $v$ of $T$ 
was always sent into $\gamma(v)$. Nevertheless, the two groups 
${\mathcal A}_T$ in \cite{KS} and the present paper coincide. Notice that there is a homeomorphism between the two differently punctured surfaces which 
slide all the punctures  of the pending lines one unit such that their 
first punctures belong are now the vertices of $T$. This homeomorphism 
conjugates between the two versions of ${\mathcal A}_T$. 
Our version has the advantage of simplifying the already cumbersome 
computations of the next section. 
\end{remark}

\subsection{The abelianized extension $\A$}
 
\begin{proposition}\label{class}
The class  
$c_{\A}\in H^2(T)$ is given by $c_{\A}=\alpha$. 
\end{proposition}
\begin{proof}
This is already stated in \cite{GrS}. In fact ${\mathcal A}_T$ splits 
over the cyclic subgroups $\Z/2^m\Z\subset T$ for all $m$ and this implies 
that the coefficient of $\chi$ vanishes. Moreover, the coefficient of 
$\alpha$ is shown in \cite{GrS} to be one.  
\end{proof}

\begin{proposition}\label{coeff}
There is an isomorphism of extensions between $\A$ and 
$T_{30, 16, 3, 1}$. 
\end{proposition} 
\begin{proof}
We have to consider $\alpha$ and $\beta$ as elements of the 
group of homeomorphisms preserving the dyadics of $S^1$. 
Recall that $T$ has the standard generators $A,B,C$ from \cite{CFP}, 
as described in section 2.1.
Realizing $A,B,C$ as (stable) couples of binary trees we can identify 
\[ \beta =C^{-1}, \,\, \alpha= C^{-1}B \]
and thus 
\[ \beta(x)=\left\{ \begin{array}{lll}
\frac{x}{2}+\frac{1}{2}, & \mbox{ if }  & x\in[0,\frac{1}{2}] \\
x+\frac{1}{4}, & \mbox{ if }  & x\in[\frac{1}{2}, \frac{3}{4}] \\
2x-\frac{3}{2}, & \mbox{ if }  & x\in[\frac{3}{4}, 1] 
\end{array}\right., \,\,\, 
 \alpha(x)=\left\{ \begin{array}{lll}
\frac{x}{2}+\frac{1}{2}, & \mbox{ if }  & x\in[0,\frac{3}{4}] \\
x+\frac{1}{8}, & \mbox{ if }  & x\in[\frac{3}{4}, \frac{7}{8}] \\
4x-\frac{7}{2}, &  \mbox{ if }  & x\in[\frac{7}{8}, 1] 
\end{array}\right. 
\] 
This implies that 
\[ \beta''(0)=-2, \beta''\left(\frac{1}{2}\right)=1, \beta''\left(\frac{3}{4}\right)=1, \,\,\, 
\alpha''(0)=-3, \alpha''\left(\frac{1}{2}\right)=0, 
\alpha''\left(\frac{3}{4}\right)=1,  
\alpha''\left(\frac{7}{8}\right)=2 \]

\vspace{0.2cm}\noindent 
In the enhanced rooted tree model we can therefore explain the 
action of $\alpha$ and $\beta$ as in the pictures below. 

\vspace{0.2cm}
\begin{center}
\includegraphics{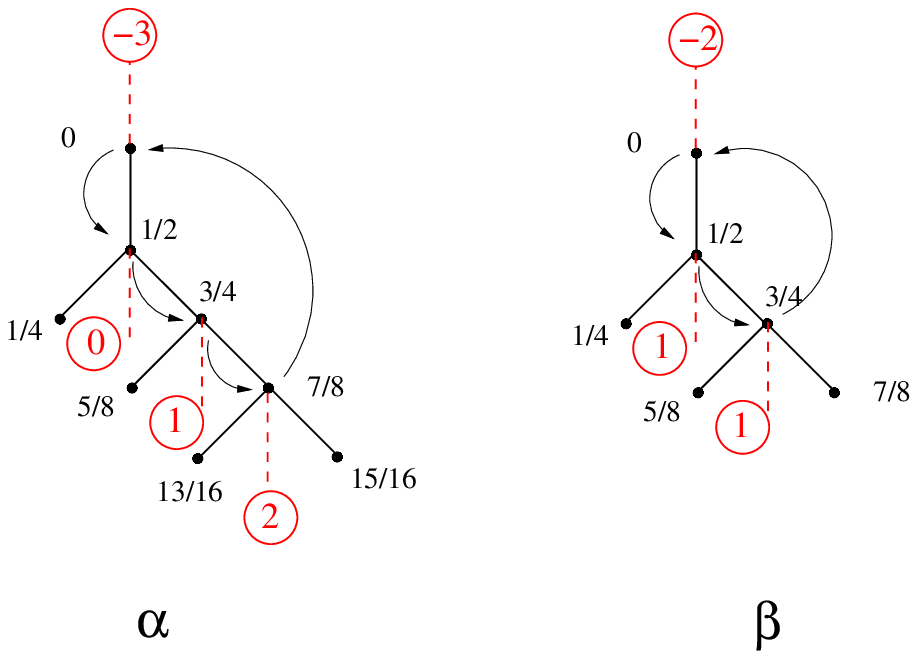}
\end{center}

\vspace{0.2cm}\noindent 
It follows by direct calculation that the action of $\beta\alpha$ 
is described by:

\vspace{0.2cm}
\begin{center}
\includegraphics{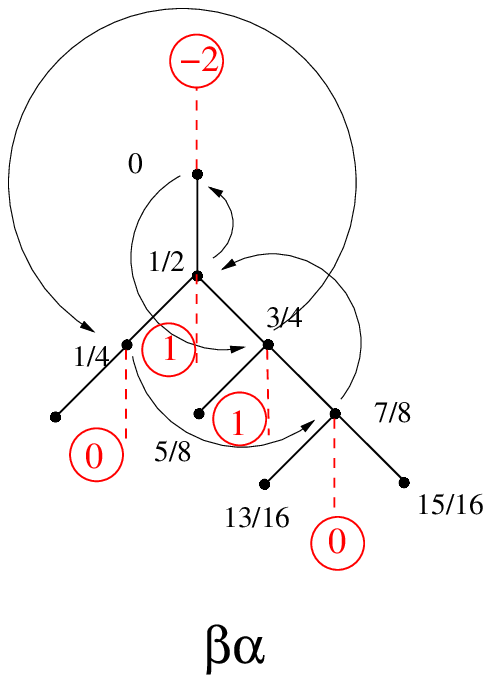}
\end{center}

\vspace{0.2cm}
\noindent
The problem we face now is to consider lifts $\alpha^{**}$ and 
$\beta^{**}$ of  
$\alpha$ and $\beta$ as mapping classes of homeomorphisms of the enhanced 
surface $ED^*$. 

\vspace{0.2cm}
\noindent
For the sake of simplicity we change in the pictures below the labels corresponding  to the vertices of 
the rooted tree, as follows: $A$ states for $0$, $B$ for $1/2$, 
$C$ for $3/4$, $D$ for $7/8$ and $E$ for $1/4$. Moreover the 
pending line at $A$ has its vertices labeled $A_1,A_2,\ldots$ and so 
on for all other vertices. 

\vspace{0.2cm}
\noindent
The supports of these classes homeomorphisms correspond to suitable disks 
around the vertices. 

\begin{enumerate}
\item $\beta^{**}$ has a support a disk embedded into $ED^*$ containing 
$A_i,B_i,C_i$ for $i\leq 3$. Moreover the action of $\beta$ is 
described as follows: first $\beta$ acts as a rotation of order 3 in the 
plane; next the vertices $A_1$ and $A_2$ are slid in counterclockwise 
direction towards the positions $C_1$ and $B_1$, respectively. 
In meantime the punctures $C_j$ (and $B_j$) are simultaneously 
translated one unit along their 
pending lines and hence $C_j$ (respectively $B_j$) will 
arrive to the position  formerly occupied by $C_{j+1}$ (respectively 
$B_{j+1}$), for all $j\geq 1$. The punctures $A_j$ (for $j\geq 3$) 
are translated simultaneously two negative units along their 
pending line and hence $A_j$ will 
arrive to the position  formerly occupied by $A_{j-2}$. 

\vspace{0.2cm}
\begin{center}
\includegraphics{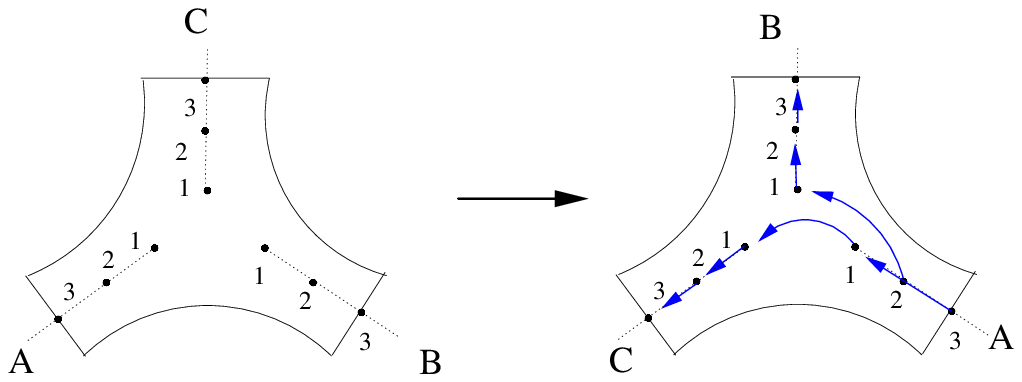}
\end{center}

\item $\alpha^{**}$  has a support a disk embedded into $ED^*$ containing 
$A_i,B_i,C_i, D_i$ for $i\leq 4$. The action of $\alpha$ is 
described as follows: first $\alpha$ acts as a rotation of order 4 in the 
plane; further $A_1$ and $A_2$ are slid in counterclockwise 
direction onto $D_1$ and respectively $D_2$, while $A_3$ 
is slid into $C_1$. The slidings occur 
simultaneously with the the translations of all punctures 
$D_j$ ($j\geq 1$) two units along their pending lines and 
the $C_j$ ($j\geq 1$) one unit along their pending line. Moreover 
the $A_j$ ($j\geq 4$) are translated three negative units along their 
pending line. The trajectories of the points are represented below 
(we did not figure the obvious translations along the pending lines):

\vspace{0.2cm}
\begin{center}
\includegraphics{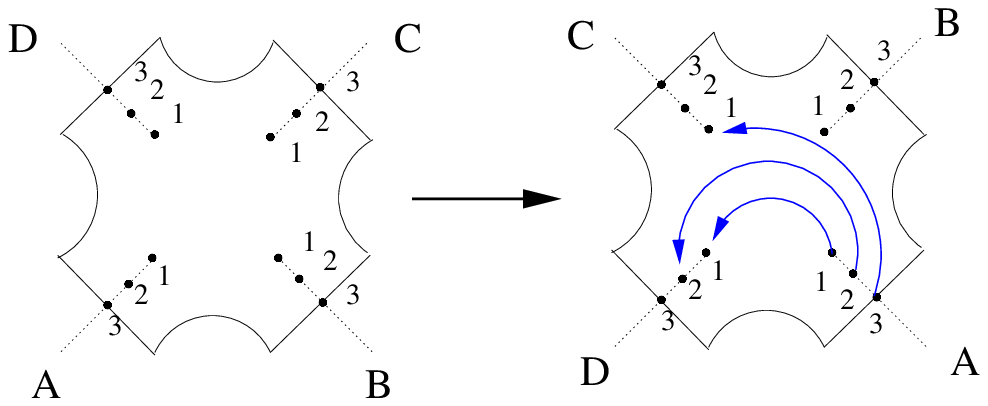}
\end{center}

\end{enumerate}

\vspace{0.2cm}
\noindent
We are able now to figure out the element 
$\beta^{**}\alpha^{**}\in {\mathcal A}_T$. Its support is now a disk 
embedded into $ED^*$ containing $A_i,B_i,C_i,D_i,E_i$ for $i\leq 4$. 

\vspace{0.2cm}
\noindent
Notice however that  the shape of the  punctures trajectories  
is not sufficient for recovering the relative position of 
punctures.  We have to specify somewhat the speed of each puncture  
along its trajectory, or equivalently, to specify a parameterization.  
There is a way to  give a discrete parameterization by associating natural numbers 
to arcs of trajectories as follows. 
The time interval is divided into 
$N$ smaller intervals for some $N$. The arc $\lambda$ is given 
the label $k\in \{1,2,\ldots, N\}$ if the respective puncture travels 
along $\lambda$ precisely in the $k$-th interval of time. 
Actually this says that whenever we have two arcs 
(disjoint or not) labeled $j$ and $k$, with $j<k$ 
the respective punctures travel first along the arc $j$ and 
next along the arc $k$.    
A mapping class group element written as a word in the $\alpha^{**}$ 
and $\beta^{**}$ will lead naturally to a discrete parameterization. 
Moreover, $N$ is chosen such that the mapping class of the 
respective homeomorphism is uniquely determined by the 
discrete parameterization. 

\vspace{0.2cm}
\noindent
As an example $\beta^{**}\alpha^{**}$ 
is completely described by the following picture 
(including parameterization):

\vspace{0.2cm}
\begin{center}
\includegraphics{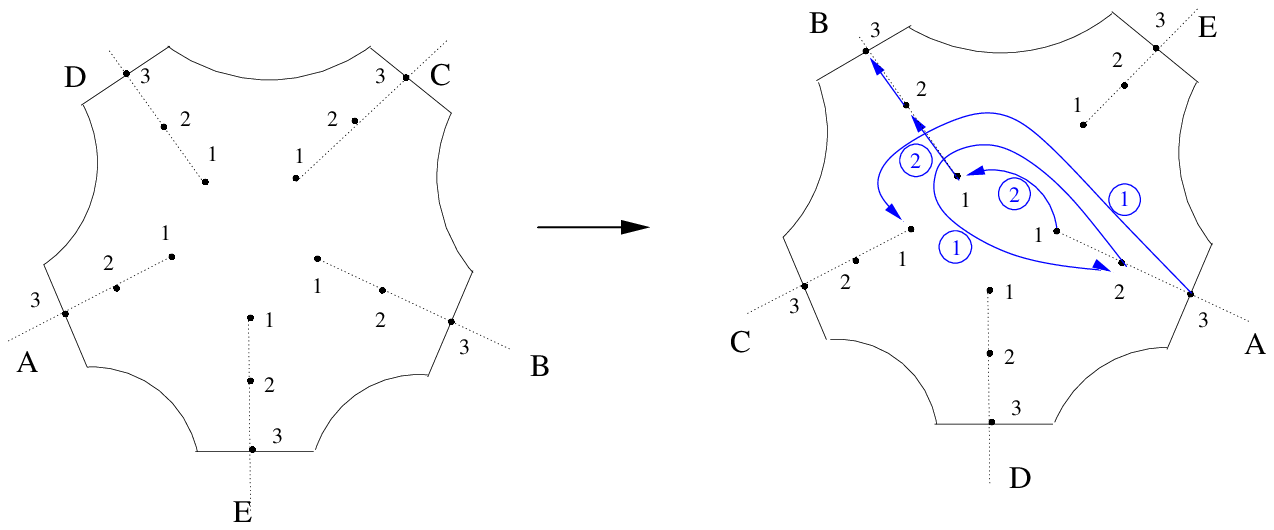}
\end{center}

\vspace{0.2cm}
\noindent
Recall that lifts in ${\mathcal A}_T$ of the relations in $T$ 
should give elements of the infinite braid group $B_{\infty}$, the  
kernel of the projection ${\mathcal A}_T\to T$. 
Denote by $a:{\mathcal A}_T\to \A$ the projection homomorphism. 
The restriction $a|_{B_{\infty}}:B_{\infty}\to Z$ is then identified  
to the abelianization homomorphism. 
We have therefore to compute the integers $n,p,q$ so that  
$a({\beta^{**}}^3)=z^q$, $a({\alpha^{**}}^4)=z^p$, 
$a((\beta^{**}\alpha^{**})^5)=z^n$.

\begin{lemma}
The braid  ${\beta^{**}}^3\in B_{\infty}$  is given by the picture 
below
  
\vspace{0.2cm}
\begin{center}
\includegraphics{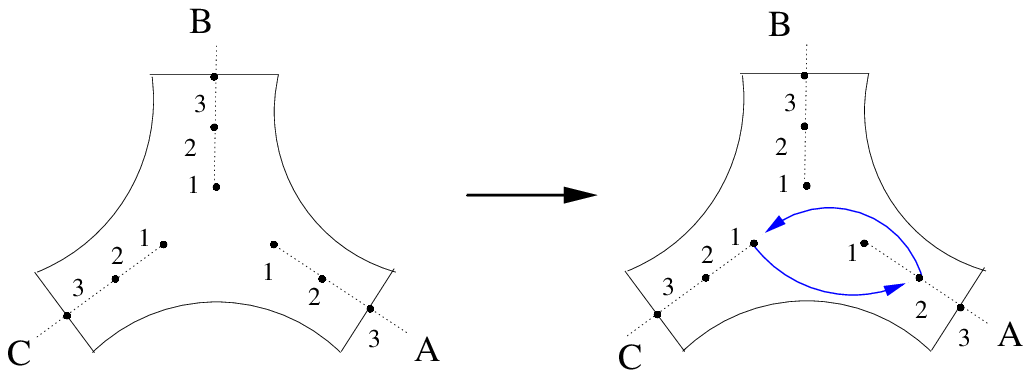}
\end{center}
In particular $a({\beta^{**}}^3)=z^3$.
\end{lemma}
\begin{proof}
We will give an explicit proof of what is going on in this (simplest) 
situation. The terminology is somewhat unconventional. We say that punctures 
``travel'' from one location to another, in certain intervals of time.
One chooses the time intervals so that their  simultaneous trajectories do not intersect and 
one  can recover the class of the associated homeomorphism.  
It is of course sufficient to have a finite number of such intervals,
which we also call steps,  and 
the relative speed within each interval is not important, because  any two 
speed values lead to isotopic homeomorphisms. Eventually we can compose 
the classes obtained by describing them step by step. 
The labels are assigned to 
the punctures and thus they travel around; at each step there are induced 
(infinite) permutations of the labels. At the end we find an element of 
a finite braid group inside $B_{\infty}$ and trajectories of punctures are 
now viewed as strands of that braid. Strands are said trivial if they are 
trivial as braid strands.  We will exemplify below with 
${\beta^{**}}^3$.

\vspace{0.2cm}\noindent 
Each time we have an action of $\beta$ there is a first step 
comprising an order 3 rotation in the 
plane and a second step, the sliding,  where  
two punctures keep traveling in 
the counterclockwise direction, while the others are fixed, 
along a circle arc of angle $\frac{4\pi}{3}$ (a $\frac{2}{3}$-turn)
and respectively $\frac{2\pi}{3}$, according to the picture of $\beta^{**}$:
the puncture labeled $A_1$ is the one which travel faster while $A_2$ 
travels a shorter amount. There are also some translations along the pending 
lines so that, for instance, $B_1$ is sent into $B_2$ and $C_1$ onto $C_2$. 
In general, 
we will not bother to represent on the picture these translations, except 
when their final action is nontrivial.
 
\vspace{0.2cm}\noindent 
Assume now that we want to compose two such classes $\beta^{**}$. 
Then, we draw first the result of the first $\beta^{**}$ action and 
further we resume with the first step, namely  an order 3 rotation.  
One applies next the second step sliding and observe 
that the puncture labeled $A_1$ is again located in the 
position from which it has to be slid.
Its fellow traveler is  this time the puncture labeled $C_1$,  
which has moved at the  previous step from position $C_1$ onto $C_2$.

\vspace{0.2cm}\noindent  
It is clear now how to add one more $\beta^{**}$: 
an  order 3 rotation sends the puncture labeled $A_1$ again 
into the position to be slid and thus it completes one more 
$\frac{2}{3}$-turn to arrive into its initial position. 
Its fellow traveler in the sliding step is the  
puncture labeled $A_2$, which followed before the path 
from $A_2$ to $B_1$ and then from $B_1$ onto the location of $B_2$. 
The sliding will send it into the location of $C_1$.

\vspace{0.2cm}
\begin{center}
\includegraphics{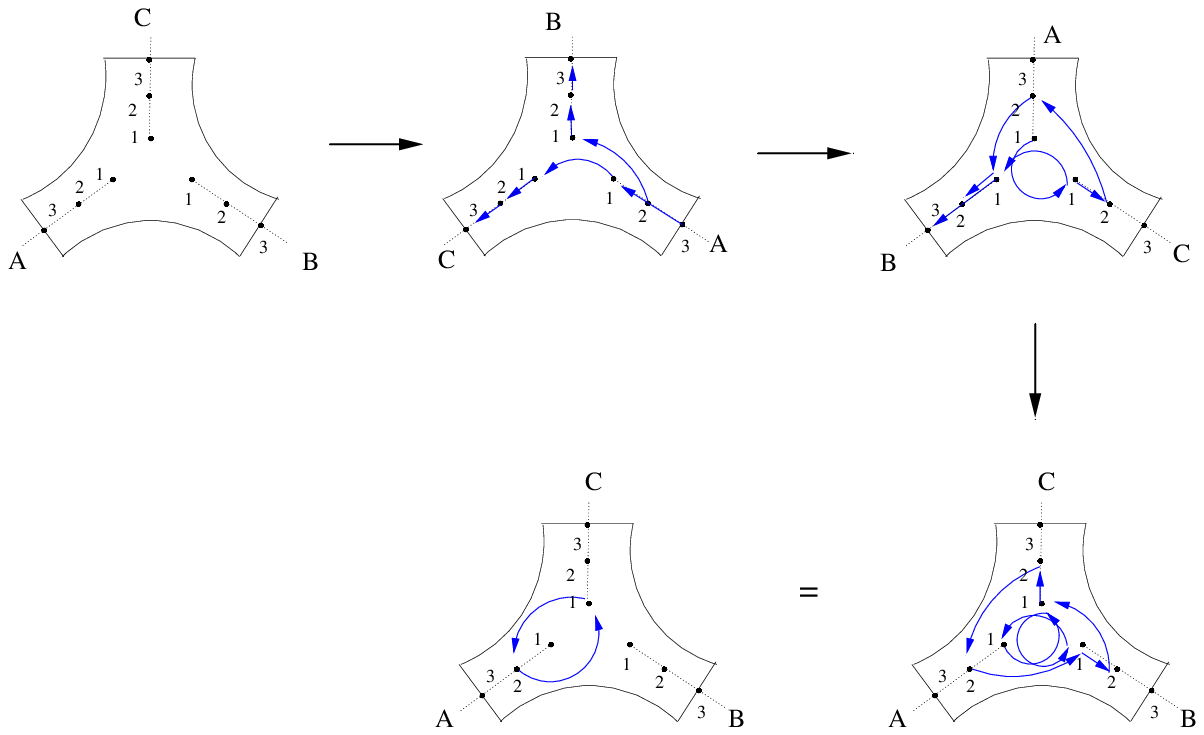}
\end{center}

\vspace{0.2cm}\noindent 
The trajectory  of $A_1$ is a strand of the braid 
${\beta^{**}}^3$ which can be split off. 
This implies that, by means of an isotopy we can assume that $A_1$ is fixed.
This isotopy corresponds to shrinking the trajectory of $A_1$ to a point. 
This shrinking could be done without touching the others trajectories, 
which means that it lifts to an isotopy between braids in the 
three dimensional space.    
It follows that the only nontrivial part of the braid ${\beta^{**}}^3$ is 
the exchange between the punctures $A_2$ and $C_1$. The translations 
along the pending lines yield trivial strands for the remaining punctures. 
\end{proof}

\end{proof}

\vspace{0.2cm}
\noindent
The remaining calculations are of the same sort, but involve more 
complicated braids. We were unable to find the braids (as the two dimensional 
picture is misleading) but we will make use of  
additional simplifications to help computing 
the images under the abelianization map.

\vspace{0.2cm}
\noindent
First, one can find $a(\sigma)$ using only  
the  winding numbers of the trajectories of the braid $\sigma$. 
Let $\sigma$  be given as a geometrical braid in $\R^2\times [0,1]$ 
by the parameterizations 
$(x_k(t),t)$, for $t\in [0,1]$, each 
subscript $k$ corresponding to one strand.
The (relative) winding number $\nu(j,k)$ of the strands $j$ 
and $k$ is  the angle that the vector 
$x_j(t)-x_k(t)$ swept when $t$ goes from $0$ to $1$. Turning 
counterclockwise yields positive angles. If the punctures sit all 
on the same line in $\R^2\times 0$ then relative winding numbers are 
multiples of $\pi$. However, they make sense even when the punctures 
are given arbitrary positions in the plane. 
The (total) winding number $\nu(\sigma)$ is the sum of all relative 
winding numbers of its (distinct) strands.

\begin{lemma}
We have $\pi a(\sigma)=\nu(\sigma)$.
\end{lemma}
\begin{proof}
Both sides are group homomorphisms and  
their values on the generators coincide. 
\end{proof}

\vspace{0.2cm}
\noindent 
It is much simpler to compute winding numbers when trajectories are 
known.  

\begin{lemma}
We have $\nu({\alpha^{**}}^4)=16 \pi$. 
\end{lemma}
\begin{proof}
We draw each trajectory individually, along with its 
discrete parameterization. The only nontrivial trajectories 
are those of the strands starting at $A_1,A_2, A_3$ and $D_1$. 
The remaining ones are easily shown to be trivial and having 
zero winding numbers with all others.

\vspace{0.2cm}
\begin{center}
\includegraphics{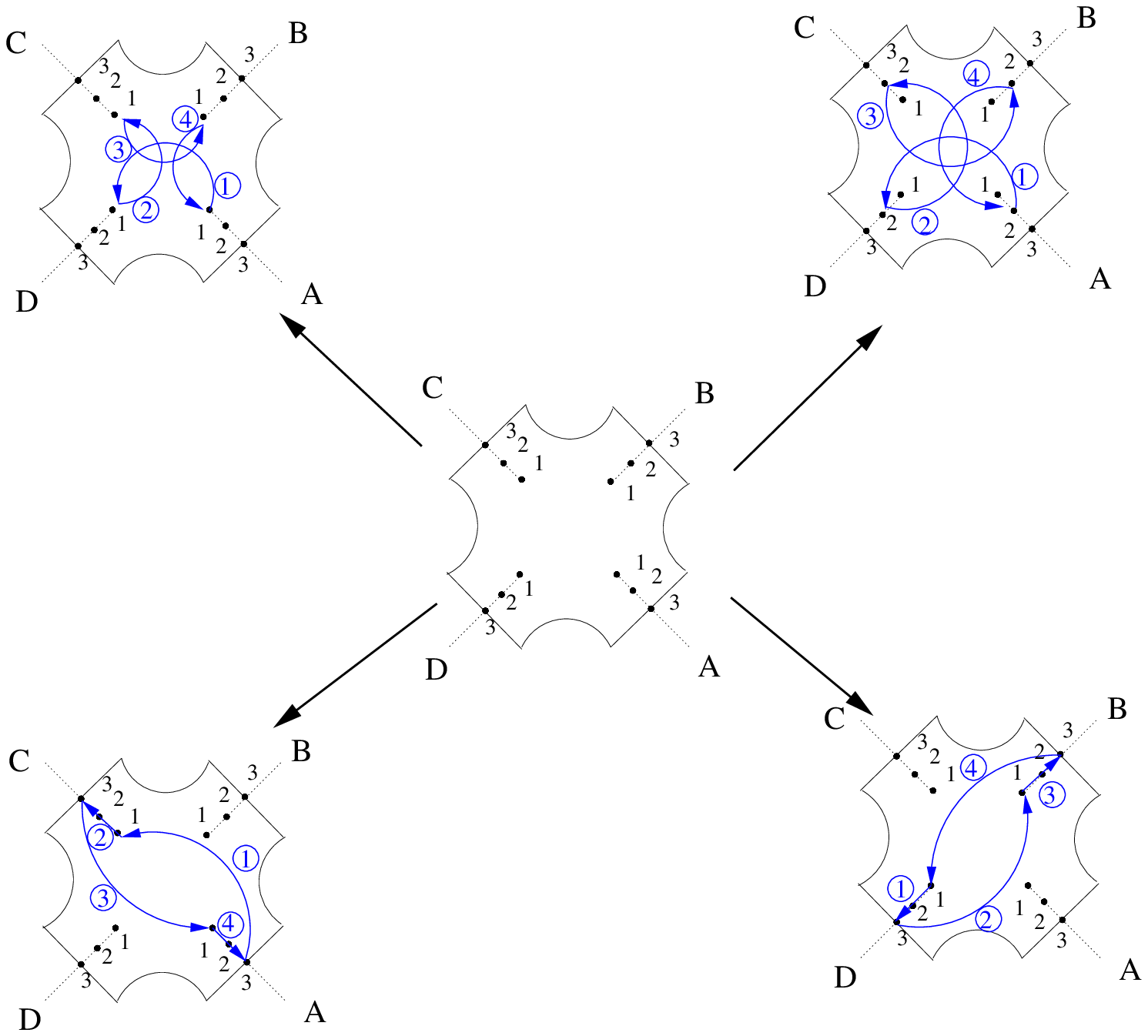}
\end{center}

\vspace{0.2cm}\noindent 
Then we compute easily from the picture above the relative 
winding numbers between strands (hereby identified to their start-point): 
\[ \nu(A_1,A_2)=6\pi, \, 
\nu(A_1,A_3)=\nu(A_2,A_3)=\nu(A_3,D_1)=\nu(A_2,D_1)=\nu(A_1,D_1)=2\pi \]
This ends the proof of the Lemma. 
\end{proof}

\begin{lemma}
We have $\nu((\beta^{**}\alpha^{**})^5)=30 \pi$. 
\end{lemma}
\begin{proof}
Recall that pictures define mapping classes of homeomorphisms.  
Then, for any mapping classes $H_1$ and $H_2$  for which the 
compositions below make sense and represent elements of $B_{\infty}$ 
we have the identity 
\[ a\left(H_1 \circ \includegraphics{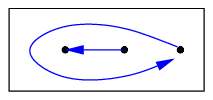}\circ H_2\right)-
a\left(H_1 \circ \includegraphics{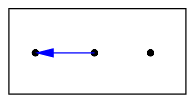}\circ H_2\right)=2 \]
In fact the first braid is obtained from the second one by inserting 
some $\sigma^2$, which braids two punctures twice. 
Replacing the first braid by the second one will be called a 
(direct) simplification of the braid diagram. 

\vspace{0.2cm}
\noindent
We can use a direct simplification within the picture 
of $\beta^{**}\alpha^{**}$ in order to remove the loop trajectory 
based at $A_2$. Thus after 5 simplifications 
within $(\beta^{**}\alpha^{**})^5$ we obtain the following trajectories.

\vspace{0.2cm}
\begin{center}
\includegraphics{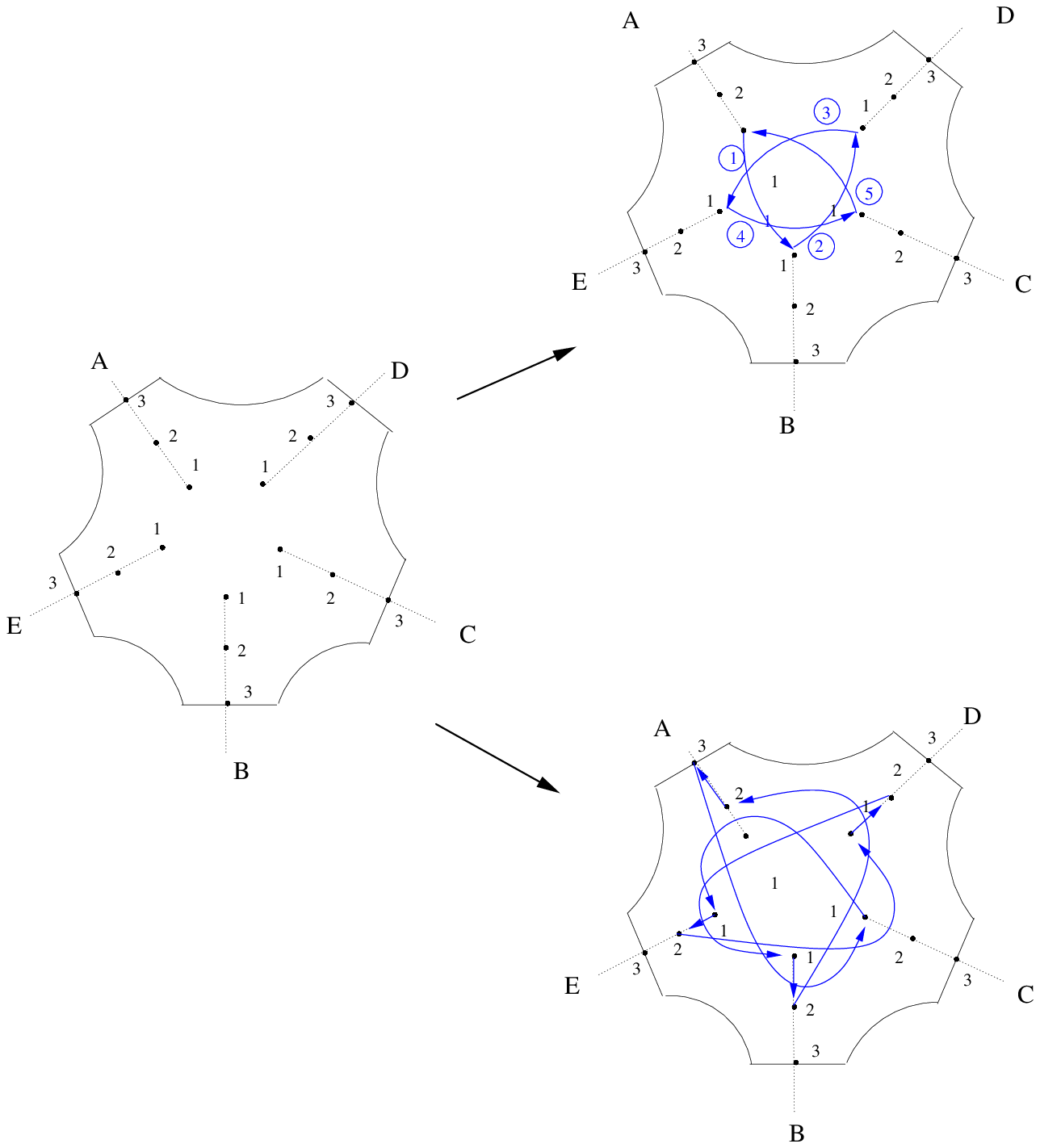}
\end{center}

\vspace{0.2cm}
\noindent
We did not draw but the nontrivial trajectories which are of two 
types: a long trajectory of the strand starting at $A_1$ and a braid cycle 
permuting circularly the nine punctures $A_2,A_3, C_1, E_1,E_2, D_1,D_2, B_1$ 
and $B_2$. The braid cycle admits 4  direct simplifications by sliding 
$C_1$ outward to $E_2D_1$, $B_1$ outward $A_3C_1$, 
$E_1$ outward of $D_2B_1$ and $D_1$ outward of $B_2A_2$ and becomes the 
following braid $b$ pictured below:

\vspace{0.2cm}
\begin{center}
\includegraphics{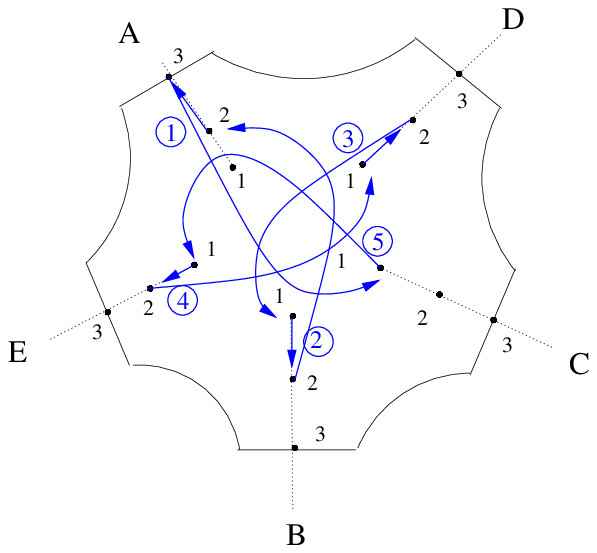}
\end{center}
  
\vspace{0.2cm}
\noindent
By direct calculation (taking care of the discrete parameterizations 
for both $b$ and the trajectory of $A_1$) we obtain  
\[ \sum_{X}\nu(A_1, X)=4\pi\]
where the sum is over all strands $X$ of the braid $b$.
Since the braid $b$ is now a cycle of length 9 we find that 
$\nu(b)=8$.  
Summing up all contributions to $\nu((\beta^{**}\alpha^{**})^5)$ 
we obtain 30, as claimed. 
\end{proof}

\begin{remark}
The elements $\beta,\alpha, \beta\alpha$ are respectively of order 
3, 4 and 5. Since the Euler class of ${\mathcal A}_T$ is zero, 
it splits over cyclic subgroups of $T$ and the same holds then for 
$\A$. This implies that $q\equiv 0({\rm mod}\, 3)$, 
$p\equiv 0({\rm mod}\, 4)$, $n\equiv 0({\rm mod}\, 5)$.
Moreover, the permutation induced by $(\beta^{**}\alpha^{**})^5$ is 
a cycle of length 9 and thus is even. In particular 
$n\equiv 0({\rm mod}\, 2)$ and so $n\equiv 0({\rm mod}\, 10)$.
\end{remark}

\vspace{0.2cm}
\noindent
\subsection{End of the proof of Theorem \ref{exte0}} 
In order to finish the proof we need to compute $\chi(n,p,q,r)$. 
We have first: 

\begin{proposition}\label{first}
We have $\chi(n,p,q,r)= 12n-15p-20q +cr$, for some $c\in \Z$. 
\end{proposition}
\begin{proof}
It suffices to use Proposition \ref{linearform} and the  
linearity of the function $\chi(n,p,q,r):\Z^4\to \Z$. For the later, use  
the arguments from the proof of lemma \ref{linearm} where the case $r=0$ is treated. 
We skip the details.
\end{proof}

\vspace{0.2cm}
\noindent 
Propositions \ref{class} and \ref{coeff} show that
$\chi(30,16,3,1)=0$, and hence $\chi(n,p,q)=12n-15p-20q-60r$.
Then proposition \ref{alpha} finishes the proof of theorem \ref{exte0}.

\vspace{0.2cm}
\noindent 
The analysis above can be used in establishing also the following: 
\begin{proposition}
The group ${\mathcal A}_{T,K}$ is finitely generated.  
\end{proposition}
\begin{proof}
Consider the group ${\mathcal A}_T$, the general case 
following the same way. 
The elements $\alpha^{**}$ and $\beta^{**}$ generate the 
quotient $T$ and it suffices to add sufficiently many elements 
to be able to generate all of $B_{\infty}$, which is the braid 
group on the punctures. Let  
$\sigma$ be the braiding of the first two punctures of 
the pending line at $A$. Let denote by $t$ the mapping class 
of the homeomorphism which translates one unit along the line made by 
gluing together the two (half)-lines pending at $A$ and $B$. 
Observe now that $t$ conjugate $\sigma$ to any of the braid generators 
of the pending line at $A$ and $B$. On the other hand 
the action of $T$ by conjugation on ${\mathcal A}_T$ is transitive 
on the set of pending lines, as it is the action of $PSL(2,\Z)$ 
on the binary tree. This means that for vertex of 
the rooted binary tree there exists a word $w=w(\alpha^{**},\beta^{**})$ 
such that the corresponding element of ${\mathcal A}_T$  
sends the pending line at $v$ asymptotically onto the pending line 
at $A$. The meaning of the word asymptotically is that all but finitely 
many vertices of the respective line are sent into the other ones. 
The first few vertices might be sent onto other pending lines, as it 
happens with $\alpha^{**}$ and $\beta^{**}$ which slide 
finitely many points. However sliding of punctures occur only 
at the supports of $\alpha^{**}$ and $\beta^{**}$. Thus the line 
pending at $v$ is sent into the line pending at $A$ and slidings of its 
vertices could appear only when its first vertex will reach 
the set $\{A_1,B_1,C_1,D_1,E_1\}$.  Let then $M$ be the mapping 
class group of an admissible subsurface containing the supports of 
$\alpha^{**}$ and $\beta^{**}$. It follows  that  
adjoining the generators of $M$ to  $\alpha^{**}, \beta^{**},\sigma, t$ 
we generate all of $B_{\infty}$ and thus ${\mathcal A}_T$. 
\end{proof}

\section{Odds and ends}
\subsection{Geometric extensions}
We would like to understand all extensions 
\[ 1 \to B_{\infty}\to G\to T\to 1 \]
coming out from {\em nature}. A tentative approach is to say that such 
an extension is geometric if there exists a tessellation of 
a planar surface $\Sigma$ with infinitely many punctures such that 
$G$ is the asymptotic mapping class group of $\Sigma$ with 
this extra structure (see also \cite{F}). Then $B_{\infty}$ is the braid group in the 
punctures. 

\vspace{0.2cm}
\noindent
In order to avoid trivial constructions we restrict to those 
examples which are {\em minimal} in some sense. The simplest 
minimality condition is to ask  
the natural homomorphism $T\to Out(B_{\infty})$ to have one orbit of 
generators of $B_{\infty}$ i.e. that $T$ acts transitively on the 
first homology of the surface $\Sigma$. Alternatively, this amounts 
to require that the lifts of the generators $\alpha$ and $\beta$ of $T$ 
together with a standard braid generator $\sigma$ form a generator 
system for $G$.

\vspace{0.2cm}
\noindent
The groups $T_{1,0,0}$ and $T_{3,1,0}$ are geometric (see \cite{FK2}) and 
minimal.  
Although  ${\mathcal A}$  is also geometric one needs to modify 
the minimality condition above in order to be fit for it.

\vspace{0.2cm}
\noindent 
It seems that there are only finitely many such minimal 
geometric extensions, for appropriate minimality conditions. 

\subsection{Finite surfaces}
It is known that the mapping class group ${\cal M}(\Sigma)$ 
of a punctured surface $\Sigma$ embeds into the groupoid of 
flips acting on the triangulations of $\Sigma$ with vertices at  
punctures.  
Quantization of the Teichm\"uller space of the surface 
$\Sigma$
lead  then by the technology of \cite{FG} 
to projective representations 
of the mapping class group ${\cal M}(\Sigma)$ and thus to 
a central extension $\widehat{{\cal M}(\Sigma)}$.

\vspace{0.2cm}
\noindent 
Recall also that  $H^2({\cal M}(\Sigma))$ is freely generated  by 
the Euler class $\chi$ together with the classes corresponding 
to each one of the punctures. 

\vspace{0.2cm}
\noindent 
It seems plausible that 
the class $c_{\widehat{{\cal M}(\Sigma)}}$ of the extension  is actually 
equal to $12 \chi \in H^2({\cal M}(\Sigma))$. 

\vspace{0.2cm}
\noindent 
Notice that  additional work is needed for obtaining this result 
because of our lack of knowledge of the 
Ptolemy groupoids and their associated groups (see \cite{pe}) for 
finite surfaces of positive genus.

{
\small      
      
\bibliographystyle{plain}

\begin{thebibliography}{30}      

\bibitem{Bar}
E.W.Barnes, {\em The genesis of the double gamma function}, Proc. London Math. Soc. 31(1899), 358-381. 


\bibitem{BB}
S.Baseilhac and R.Benedetti, {\em Classical and quantum dilogarithmic invariants of flat PSL(2,C)-bundles over 3-manifolds}, Geometry $\&$  Topology  9(2005), 493-569. 

\bibitem{BBL}
Hua Bai, F.Bonahon and Xiaobo  Liu, {\em  Local representations of the 
quantum Teichm\"uller space}, 
math.GT/0707.2151.  


\bibitem{Bax}
R.Baxter, Exactly solvable modles in statistical mechanics, Academic Press, 1982. 

\bibitem{BL}
F.Bonahon and Xiaobo Liu, {\em  Representations of the quantum 
Teichm\"uller space and invariants of surfaces diffeomorphisms}, 
Geometry $\&$  Topology  11(2007), 889--937.


\bibitem{Bri0}
M.G.Brin, {\em  The chameleon groups of Richard J. Thompson: 
automorphisms and dynamics}, 
Inst. Hautes \'Etudes Sci. Publ. Math. No. 84, 1996, 5-33.

\bibitem{Bri}
M.G.Brin, {\em The Algebra of Strand Splitting. I. A Braided Version of Thompson's Group V}, 
 J. Group Theory  10  (2007),  no. 6, 757--788. 

\bibitem{Bri2}
M.G.Brin, {\em The Algebra of Strand Splitting.II. A Presentation for the Braid Group on One Strand}, Internat. J. Algebra and Comput. 16(2006), 203-219. 




\bibitem{Ca}
D. Calegari, {\em Circular groups, Planar groups and the Euler class}, 
Proceedings of the Casson Fest, Geom. Topol. Monogr. 7(2004), 431-491. 

      
\bibitem{CFP}      
J.W.Cannon, W.J. Floyd, and W.R. Parry,       
{\em Introductory notes on Richard Thompson's groups},       
Enseign. Math. 42(1996), 215-256.  
 
\bibitem{CF}
L.Chekov and V.Fock, {\em Quantum Teichm\"uller space}, 
math.QA/9908165. 



\bibitem{De1}
P.Dehornoy, {\em Geometric presentations for Thompson's groups}, 
 Journal Pure Appl. Algebra, 203(2005), 1-44. 


\bibitem{De2}
P.Dehornoy, {\em The group of parenthesized braids},  
Advances Math. 205(2006), 354-409.


\bibitem{Fad}
L.D.Faddeev, {\em Discrete Heisenberg-Weyl group and modular group}, 
Lett. Math. Phys. 34(1995), 249-254. 
      

\bibitem{FaK}
L.Faddeev and R.Kashaev, {\em Quantum dilogarithm}, 
Mod. Phys. Lett. A  9(1994), 427--434. 



\bibitem{Fo}
V.Fock, {\em  Dual Teichm\"uller spaces}, math.DG-GA/9702018.

\bibitem{FG1}
V.Fock and A.B.Goncharov, {\em Moduli spaces of local systems 
and higher Teichm\"uller theory}, Inst. Hautes \'Etudes Sci. Publ. Math. no.103, 2006, 1-211. 

 \bibitem{FG2}
V.Fock and A.B.Goncharov, {\em Moduli spaces of convex projective structures 
on surfaces}, Advances Math.   208(2007),  249--273. 
 
\bibitem{FG}
V.Fock and A.B.Goncharov, {\em The quantum dilogarithm and unitary 
representations of cluster modular groupoids}, 
Invent. Math.  175(2009),   223--286. 

\bibitem{FG3}
V.Fock and A.B.Goncharov, {\em Cluster ensembles, quantization and the dilogarithm II: The intertwinner}, arXiv:math/0702398. 


\bibitem{Fun}
L.Funar, {\em Ptolemy groupoids actions on Teichmuller spaces}, 
Modern Trends in Geometry and Topology, Deva, Romania 2005, (D.Andrica, P.Blaga, 
S.Moroianu, Eds.), Cluj Univ.Press 2006, p.185-201.


\bibitem{F}
L.Funar, {\em Braided Houghton groups as mapping class groups},  
Annales Sci. Univ. "A.I.Cuza" Jassy, special vol. to the memory of 
Gh.Ionesei,  53(2007), 229-240.

\bibitem{FK2}
L.Funar and C.Kapoudjian, {\em The braided Ptolemy-Thompson group is finitely presented}, 
Geometry $\&$ Topology 12(2008), 375-430. 

\bibitem{FK3}
L.Funar and C.Kapoudjian, {\em The Ptolemy-Thompson group $T^*$ is 
asynchronously combable}, math.GT/0602490.






    
\bibitem{GS}      
E.Ghys and V.Sergiescu, {\em Sur un groupe remarquable de      
  diff\'eomorphismes du cercle}, Comment.Math.Helv. 62(1987), 185-239.      
 


\bibitem{Go}
A.B.Goncharov, {\em Pentagon relation for the quantum dilogarithm and quantized 
${\cal M}\sb {0,5}\sp {\rm cyc}$},   
Geometry and dynamics of groups and spaces, 
In memory of Alexander Reznikov (M.Kapranov, S.Kolyada, Yu.I.Manin, P.Moree and L.Potyagailo Ed.), 415--428, Progr. Math., 265, Birkh\"auser, 
Basel, 2008.

\bibitem{GrS}
P.Greenberg and V.Sergiescu, {\em 
An acyclic extension of the braid group}, 
Comment. Math. Helv.  66(1991),  109-138. 


\bibitem{K}
R.Kashaev, {\em Quantization of Teichm\"uller spaces and quantum 
dilogarithm}, Lett. Math. Phys. 43(1998), 105-115. 


\bibitem{KS}
C.Kapoudjian and V.Sergiescu, {\em An extension of the Burau representation to a mapping class group associated to Thompson's group $T$},
Geometry and dynamics,  141--164, Contemp. Math., 389, Amer. Math. Soc., Providence, RI, 2005.   
      

     
\bibitem{LS}      
P. Lochak and L. Schneps,      
\newblock {\em The universal Ptolemy-Teichm\"uller groupoid},       
\newblock { in  Geometric Galois actions}, vol. 2,  L.M.S. Lecture      
Notes Ser., 243, Cambridge Univ.Press, 1997.       
      
 
 

  
     
\bibitem{pe0}      
R.C.Penner, {\em Universal constructions in Teichmuller theory}, 
Advances  Math.   98(1993),  143-215.      
      
\bibitem{pe}       
R.C.Penner, {\em The universal Ptolemy group and its completions},   
Geometric Galois actions, 2, 293-312, L.M.S. Lecture Notes Ser., 243,
Cambridge Univ. Press, Cambridge, 1997.      





      
\end{thebibliography}

}

\end{document}